\documentclass[12pt]{amsart}

\usepackage{amssymb}
\usepackage{bm}
\usepackage{graphicx}
\usepackage{ifthen}
\usepackage{pict2e}
\usepackage{xargs}
\usepackage{xspace}

\newcommand{\definedterm}[1]{\emph{#1}}

\newcommand{\AD}{\mathtt{AD}}

\newcommand{\Bairespace}[1][]{
  \ifthenelse{\equal{#1}{}}{\functions{\N}{\N}}{\functions{#1}{\N}}
}
\newcommand{\Bairetree}[1][]{
  \ifthenelse{\equal{#1}{}}{\functions{<\N}{\N}}{\functions{#1}{\N}}
}
\newcommand{\ball}[3][]{\calB_{#1}(#2, #3)}
\newcommand{\bbP}{\mathbb{P}}
\newcommand{\bbQ}{\mathbb{Q}}
\newcommand{\bbS}{\mathbb{S}}

\newcommand{\bfb}{\mathbf{b}}
\newcommand{\bfc}{\mathbf{c}}
\newcommand{\bfC}{\mathbf{C}}
\newcommand{\bfd}{\mathbf{d}}
\newcommand{\bfD}{\mathbf{D}}
\newcommand{\bfone}{\mathbf{1}}

\newcommand{\boundingnumber}{\mathfrak{b}}

\newcommand{\calB}{\mathcal{B}}
\newcommand{\calF}{\mathcal{F}}

\newcommand{\calI}{\mathcal{I}}
\newcommand{\calJ}{\mathcal{J}}
\newcommand{\calM}{\mathcal{M}}
\newcommand{\calN}{\mathcal{N}}

\newcommand{\calU}{\mathcal{U}}
\newcommand{\calV}{\mathcal{V}}

\newcommand{\Cantorspace}[1][]{
  \ifthenelse{\equal{#1}{}}{\functions{\N}{2}}{\functions{#1}{2}}
}
\newcommand{\Cantortree}[1][]{
  \ifthenelse{\equal{#1}{}}{\functions{<\N}{2}}{\functions{#1}{2}}
}
\newcommand{\cardinality}[1]{|#1|}
\newcommand{\CH}{\mathtt{CH}}
\newcommand{\characteristicfunction}[1]{\bfone_{#1}}

\newcommand{\chromaticnumber}[1]{\chi(#1)}
\newcommand{\closedinterval}[2]{[#1, #2]}
\newcommand{\closure}[1]{\overline{#1}}
\newcommand{\cofinality}[1]{\mathrm{cof}(#1)}
\newcommand{\compactification}[1]{#1_*}

\newcommand{\completegraph}[1]{\mathbb{K}_{#1}}

\newcommand{\composition}{\circ}

\newcommandx{\concatenation}[2][1 = undefined, 2 = undefined]{
  \ifthenelse{\equal{#1}{undefined}}{{}\smallfrown}{
    \ifthenelse{\equal{#2}{undefined}}{\bigoplus #1}{\bigoplus_{#1} #2}
  }
}
\newcommandx{\constant}[2][2 =]{\ifthenelse{\equal{#2}{}}{c_{#1}}{c_{#1, #2}}}
\newcommandx{\constantfunction}[3][2 =, 3 =]{
  \ifthenelse{\equal{#2}{}}{c \from #1 \to \image{c}{#1}}{c_{#3} \from #1 \to #2}
}
\newcommand{\constantsequence}[2]{\sequence{#1}^{#2}}
\newcommand{\continuum}{\mathfrak{c}}

\newcommandx{\convolution}[2][1 = undefined, 2 = undefined]{
  \ifthenelse{\equal{#1}{undefined}}{\mathrel{*}}{
    \ifthenelse{\equal{#2}{undefined}}{\bigotimes #1}{\bigotimes_{#1} #2}
  }
}
\newcommand{\coveringnumber}[1]{\mathrm{cov}(#1)}
\newcommand{\DC}{\mathtt{DC}}
\newcommandx{\Deltaclass}[2][1=,2=]{
  \ifthenelse{\equal{#2}{}}{\mathbf{\Delta}_{#1}}{\mathbf{\Delta}^{#1}_{#2}}
}
\newcommand{\diagonal}[2][]{\triangle^{#1}(#2)}

\newcommandx{\directsum}[2][1 = undefined, 2 = undefined]{
  \ifthenelse{\equal{#1}{undefined}}{{}\oplus}{
    \ifthenelse{\equal{#2}{undefined}}{\bigoplus #1}{\bigoplus_{#1} #2}
  }
}

\newcommandx{\disjointunion}[2][1 =, 2 =]{
  \ifthenelse{\equal{#1}{}}{\sqcup}{
    \ifthenelse{\equal{#2}{}}{\bigsqcup #1}{{\bigsqcup_{#1} #2}}
  }
}
\newcommand{\doesnotgoto}{\not\goesto}
\newcommand{\domain}[1]{\mathrm{dom}(#1)}
\newcommand{\dominatingnumber}{\mathfrak{d}}

\newcommand{\emptysequence}{\emptyset}
\newcommand{\eventuallydominatedby}{\le^*}

\newcommand{\extendedby}{\sqsubseteq}

\newcommand{\extension}[2][]{\ifthenelse{\equal{#1}{}}{\overline{#2}}{\overline{#2}_{#1}}}
\newcommandx{\extensions}[2][2 =]{
  \ifthenelse{\equal{#2}{}}{\calN_{#1}}{\calN_{#1, #2}}
}

\newcommand{\forces}[1][]{\Vdash_{#1}}
\newcommand{\forcofinitelymany}{\forall^\infty}
\newcommand{\from}{\colon}
\newcommand{\Fsigma}{$F_\sigma$\xspace}
\newcommandx{\functions}[3][3 =]{
  \ifthenelse{\equal{#3}{}}{#2^{#1}}{#2^{#1}_{#3}}
}

\newcommand{\Gdelta}{$G_\delta$\xspace}
\newcommand{\goesto}{\rightarrow}
\newcommand{\graph}[1]{\mathrm{graph}(#1)}

\newcommand{\GzeroN}[1][]{
  \ifthenelse{\equal{#1}{}}{\mathbb{G}_0^\N}{\mathbb{G}_{0, #1}^{\N}}
}
\newcommand{\Hcountable}[1][]{\mathbb{H}_{\functions{\N}{(#1 \times \N)}, \infty}}

\newcommand{\Heven}[1][]{\mathbb{H}_{\functions{\N}{#1}}}
\newcommand{\Hevenprime}[1][]{\mathbb{H}_{\functions{\N}{#1}}'}

\newcommand{\image}[2]{#1(#2)}

\newcommandx{\identityfunction}[2][2 =]{
  \ifthenelse{\equal{#2}{}}{\mathrm{id} \from #1 \to #1}{\mathrm{id} \from #1 \to #2}
}
\newcommand{\infimum}[2][]{
  \ifthenelse{\equal{#1}{}}{\inf #1}{\inf_{#1}{#2}}
}
\newcommand{\injections}[1]{\functions{\N}{(\N)}}
\newcommand{\injectivefunctions}[2]{\functions{#1}{(#2)}}
\newcommandx{\intersection}[2][1 =, 2 =]{
  \ifthenelse{\equal{#1}{}}{\cap}{
    \ifthenelse{\equal{#2}{}}{\bigcap #1}{{\bigcap_{#1} #2}}
  }
}

\newcommand{\Ksigma}{$K_\sigma$\xspace}
\newcommand{\length}[1]{|#1|}
\newcommandx{\limit}[2][1 =, 2 =]{
  \ifthenelse{\equal{#1}{}}{\lim}{
    \ifthenelse{\equal{#2}{}}{\lim #1}{{\lim_{#1} #2}}
  }
}

\newcommand{\mathand}{\text{ and }}

\newcommand{\meagerideal}{\calM}

\newcommand{\modifiedpower}[3]{\mathrm{Ext}_{#3}(#2^{#1})}
\newcommand{\N}{\mathbb{N}}

\newcommand{\neventuallydominatedby}{\not\eventuallydominatedby}
\newcommand{\nforces}[1][]{\not\forces[#1]}

\newcommand{\nullideal}{\calN}
\newcommand{\ogd}[2][]{\mathtt{OGD}^{#1}(#2)}

\newcommand{\pair}[2]{(#1, #2)}

\newcommand{\partialto}{\rightharpoonup}
\newcommandx{\Piclass}[2][1=,2=]{
  \ifthenelse{\equal{#2}{}}{\mathbf{\Pi}_{#1}}{\mathbf{\Pi}^{#1}_{#2}}
}

\newcommand{\predecessor}[1]{#1^{\mathord{-}}}
\newcommand{\preimage}[2]{#1^{-1}(#2)}
\newcommandx{\product}[2][1 =, 2 =]{
  \ifthenelse{\equal{#1}{}}{\times}{
    \ifthenelse{\equal{#2}{}}{\prod #1}{{\prod_{#1} #2}}
  }
}
\newcommandx{\projection}[2][1 =, 2 =]{
  \ifthenelse{\equal{#1}{}}{\mathrm{proj}}{
    \ifthenelse{\equal{#2}{}}{\projection_{#1}}{
      \image{\projection[#1]}{#2}
    }
  }
}
\newcommand{\R}{\mathbb{R}}
\newcommand{\Ramseyspace}{\functions{\N}{[\N]}}
\newcommand{\Rcon}[2][]{
  \ifthenelse{\equal{#1}{}}
    {\mathrm{Cnvg}(#2)}
    {\mathrm{Cnvg}_{#2}(\functions{\N}{#1})}
}
\renewcommandx{\restriction}[3][3 = undefined]{
  \ifthenelse{\equal{#3}{undefined}}{#1 \upharpoonright #2}{#1 \upharpoonright #2 \to #3}
}

\newcommandx{\sequence}[2][2 = undefined]{
  \ifthenelse{\equal{#2}{undefined}}{(#1)}{
    (#1)_{#2}
  }
}
\newcommandx{\set}[2][2 = undefined]{
  \ifthenelse{\equal{#2}{undefined}}{\{ #1 \}}{
    \{ #1 \suchthat #2 \}
  }
}
\newcommand{\setcomplement}[1]{\twiddle #1}
\newcommandx{\sets}[3][3 =]{
  \ifthenelse{\equal{#3}{}}{[#2]^{#1}}{[#2]^{#1}_{#3}}
}

\newcommandx{\Sigmaclass}[2][1=,2=]{
  \ifthenelse{\equal{#2}{}}{\mathbf{\Sigma}_{#1}}{\mathbf{\Sigma}^{#1}_{#2}}
}
\newcommand{\strictlyextendedby}{\sqsubset}

\newcommand{\suchthat}{\mid}
\newcommand{\support}[1]{\text{supp} \ #1}
\newcommand{\textexponent}[2]{$\text{#1}^{\text{#2}}$}

\newcommand{\triple}[3]{(#1, #2, #3)}
\newcommand{\twiddle}{\raisebox{1pt}{\scalebox{.75}{$\mathord{\sim}$}}}
\newcommandx{\union}[2][1 =, 2 =]{
  \ifthenelse{\equal{#1}{}}{\cup}{
    \ifthenelse{\equal{#2}{}}{\bigcup #1}{{\bigcup_{#1} #2}}
  }
}

\newcommand{\ZF}{\mathtt{ZF}}
\newcommand{\ZFC}{\mathtt{ZFC}}

\newcommand{\Baire}{Baire\xspace}
\newcommand{\Borel}{Bor\-el\xspace}

\newcommand{\Cauchy}{Cau\-chy\xspace}

\newcommand{\Hausdorff}{Haus\-dorff\xspace}
\newcommand{\Hurewicz}{Hur\-e\-wicz\xspace}

\newcommand{\Jayne}{Jayne\xspace}
\newcommand{\Kechris}{Kech\-ris\xspace}
\newcommand{\Lebesgue}{Leb\-esgue\xspace}
\newcommand{\Lecomte}{Lec\-omte\xspace}
\newcommand{\Louveau}{Lou\-veau\xspace}

\newcommand{\Polish}{Po\-lish\xspace}
\newcommand{\Rogers}{Rog\-ers\xspace}
\newcommand{\Sacks}{Sacks\xspace}
\newcommand{\SaintRaymond}{Saint Ray\-mond\xspace}
\newcommand{\Silver}{Sil\-ver\xspace}
\newcommand{\Solecki}{Sol\-eck\-i\xspace}
\newcommand{\Souslin}{Sous\-lin\xspace}

\newcommand{\Woodin}{Wood\-in\xspace}
\newcommand{\Zeleny}{Zel\-e\-ny\xspace}

\newenvironment{lemmaproof}{
  
  \begin{proof}
}{\end{proof}}

\newenvironment{propositionproof}{
  
  \begin{proof}
}{\end{proof}}

\newenvironment{theoremproof}{
  
  \begin{proof}
}{\end{proof}}

\newtheorem{lemma}{Lemma}[section]
\newtheorem{proposition}[lemma]{Proposition}

\newtheorem{theorem}[lemma]{Theorem}

\theoremstyle{definition}

\begin{document}


\begin{abstract}
  We show that several dichotomy theorems concerning the second
  level of the \Borel hierarchy are special cases of the
  $\aleph_0$-dimen\-sional generalization of the open graph dichotomy,
  which itself follows from the usual proof(s) of the perfect set theorem.
  Under the axiom of determinacy, we obtain the generalizations of
  these results from analytic metric spaces to separable metric spaces.
  We also consider connections between cardinal invariants and the
  chromatic numbers of the corresponding dihypergraphs.
\end{abstract}

\author[R. Carroy]{Rapha\"{e}l Carroy}

\address{
  Rapha\"{e}l Carroy \\
  Kurt G\"{o}del Research Center for Mathematical Logic \\
  Universit\"{a}t Wien \\
  W\"{a}hringer Stra{\ss}e 25 \\
  1090 Wien \\
  Austria
}

\email{raphael.carroy@univie.ac.at}

\urladdr{
  http://www.logique.jussieu.fr/~carroy/indexeng.html
}

\author[B.D. Miller]{Benjamin D. Miller}

\address{
  Benjamin D. Miller \\
  Kurt G\"{o}del Research Center for Mathematical Logic \\
  Universit\"{a}t Wien \\
  W\"{a}hringer Stra{\ss}e 25 \\
  1090 Wien \\
  Austria
 }

\email{benjamin.miller@univie.ac.at}

\urladdr{
  http://www.logic.univie.ac.at/benjamin.miller
}

\author[D.T. Soukup]{D\'{a}niel T. Soukup}

\address{
  D\'{a}niel T. Soukup \\
  Kurt G\"{o}del Research Center for Mathematical Logic \\
  Universit\"{a}t Wien \\
  W\"{a}hringer Stra{\ss}e 25 \\
  1090 Wien \\
  Austria
 }

\email{daniel.soukup@univie.ac.at}

\urladdr{
  http://www.logic.univie.ac.at/$\twiddle$soukupd73
}

\thanks{The authors were supported in part by FWF Grants P28153,
  P29999, and I1921.}
  
\keywords{Coloring, dichotomy, separation, sigma-continuous}

\subjclass[2010]{Primary 03E15, 26A21, 28A05, 54H05}

\title
  [The open dihypergraph dichotomy]
  {The open dihypergraph dichotomy and the second level of the Borel
    hierarchy}

\maketitle

\section*{Introduction}

\subsection*{Basic notions}
A topological space $X$ is \definedterm{\Polish} if it is second
countable and completely metrizable. A topological space $X$ is
\definedterm{analytic} if it is a continuous image of a closed subset
of $\Bairespace$.

A subset of a topological space is \definedterm{\Ksigma} if it is a
countable union of compact sets, {\Fsigma} if it is a countable union of
closed sets, \definedterm{\Gdelta} if it is a countable intersection of
open sets, \definedterm{$\Deltaclass[0][2]$} if it is both \Fsigma and
\Gdelta, and \definedterm{\Borel} if it is in the smallest $\sigma$-algebra
containing the open sets.

A function $\pi \from X \to Y$ is \definedterm{$\Gamma$-measurable}
if the pre-image of every open set is in $\Gamma$. When $\Gamma$ is
the family of \Borel sets, we say that such a function is \definedterm
{\Borel}. We say that a function $\pi \from X \to Y$ is \definedterm
{$\sigma$-continuous with $\Gamma$ witnesses} if $X$ is the union
of countably-many sets in $\Gamma$ on which $\phi$ is continuous.

Given a set $C \subseteq X$, we use $\setcomplement{C}$ to denote
the complement $C$ within $X$. We say that $C$ \definedterm
{separates} a set $A \subseteq X$ from a set $B \subseteq X$ if $A
\subseteq C$ and $B \subseteq \setcomplement{C}$. A subset of
$\product[d \in D][X_d]$ is a \definedterm{hyperrectangle} if it is of the
form $\product[d \in D][Y_d]$, where $Y_d \subseteq X_d$ for all $d
\in D$. When $D = 2$, we say that such a set is a \definedterm
{rectangle}. The \definedterm{box topology} on a product $\product[d
\in D][X_d]$ of topological spaces is the topology generated by the
sets of the form $\product[d \in D][U_d]$, where $U_d \subseteq X_d$
is open for all $d \in D$.

A \definedterm{$D$-ary relation} on $X$ is a subset of $\functions
{D}{X}$. A \definedterm{homomorphism} from a $D$-ary relation
$R$ on $X$ to a $D$-ary relation $S$ on $Y$ is a map $\pi \from X
\to Y$ such that $\image{\pi^D}{R} \subseteq S$.
A \definedterm{$D$-dimensional dihypergraph} on $X$ is a $D$-ary
relation $H$ on $X$ disjoint from the set $\diagonal[D]{X}$ of constant
sequences. When $D = 2$, we say that a symmetric such set is a
\definedterm{graph}. The \definedterm{complete $D$-dimensional
dihypergraph} on $X$ is the complement of $\diagonal[D]{X}$. A
set $Y \subseteq X$ is \definedterm{$H$-independent} if $\restriction
{H}{Y} = \emptysequence$. A \definedterm{$\kappa$-coloring} of a
$D$-dimensional dihypergraph $H$ is a homomorphism from $H$ to
the complete $D$-dimensional dihypergraph on a set of cardinality
$\kappa$. The \definedterm{chromatic number} of $H$, or
$\chromaticnumber{H}$, is the least cardinal $\kappa$ for which there
is a $\kappa$-coloring of $H$. We say that an $\N$-dimensional
dihypergraph is \definedterm{hereditary} if it is closed under
subsequences.

For each partial function $t \from \N \partialto D$, set $\extensions{t}
= \set{\bfd \in \functions{\N}{D}}[t \extendedby \bfd]$. We use $\sequence
{d}$, $\constantsequence{d}{n}$, and $\constantsequence{d}{\infty}$ to
denote the sequences with constant value $d$ in $\functions{1}{D}$,
$\functions{n}{D}$, and $\functions{\N}{D}$. We use $s \concatenation t$
and $\concatenation[n \in N][s_n]$ to denote concatenation.

We use $\completegraph{\Cantorspace}$ to denote the complete
graph on $\Cantorspace$. More generally, we use $\Heven[D]$ to
denote the $D$-dimensional dihypergraph on $\functions{\N}{D}$
given by $\Heven[D] = \union[t \in \functions{<\N}{D}][{\product[d \in
D][\extensions{t \concatenation \sequence{d}}]}]$, and $\Hevenprime
[D]$ to denote the $D$-dimensional dihypergraph on $\functions{\N}{D}$
given by $\Hevenprime[D] = \union[\pair{\bfd}{t} \in
\injectivefunctions{D}{D} \times \functions{<\N}{D}][{\product[d \in D]
[\extensions{t \concatenation \sequence{\bfd(d)}}]}]$, where
$\injectivefunctions{N}{D}$ denotes the set of injective elements
of $\functions{N}{D}$.

The \definedterm{covering number} of an ideal $\calI$ on a set
$X$ is the least cardinal $\coveringnumber{\calI}$ for which $X$ is the
union of $\coveringnumber{\calI}$-many sets in $\calI$. The
\definedterm{cofinality} of an ideal $\calI$ on a set $X$ is the
least cardinal $\cofinality{\calI}$ for which there is a set $\calJ
\subseteq \calI$ of cardinality $\cofinality{\calI}$ such that
$\forall I \in \calI \exists J \in \calJ \ I \subseteq J$. We use
$\meagerideal$ to denote the $\sigma$-ideal on $\R$ consisting
of the meager sets, and $\nullideal$ to denote the $\sigma$-ideal on
$\closedinterval{0}{1}$ consisting of the \Lebesgue null sets.

We write $\forcofinitelymany n \in N \ \phi(n)$ to indicate that $\set{n
\in N}[\neg \phi(n)]$ is finite. The order of \definedterm{eventual
domination} is the quasi-order on $\Bairespace$ given by $\bfc
\eventuallydominatedby \bfd \iff \forcofinitelymany n \in \N \ \bfc(n) \le
\bfd(n)$. The \definedterm{bounding number} is the least cardinal
$\boundingnumber$ for which there is a set $\calF \subseteq
\Bairespace$ of cardinality $\boundingnumber$ such that $\forall \bfc
\in \Bairespace \exists \bfd \in \calF \ \bfd \not \eventuallydominatedby \bfc$.
The \definedterm{dominating number} is the least cardinal
$\dominatingnumber$ for which there is a set $\calF \subseteq
\Bairespace$ of cardinality $\dominatingnumber$ such that $\forall \bfc
\in \Bairespace \exists \bfd \in \calF \ \bfc \eventuallydominatedby \bfd$.

\subsection*{Motivation}
In \cite{Hurewicz}, the existence of a minimal \Borel subset of a \Polish
space that is not \Fsigma was established. In \cite{Wadge}, this was
generalized to every level of the \Borel hierarchy under \Borel
determinacy. In \cite{LouveauSaintRaymond}, these results were
established in a far weaker subsystem of $\ZF$, leading to
parametrized generalizations. In \cite{Lecomte}, they were further
generalized to subsets of products.

In \cite{LecomteZeleny}, the question as to the circumstances under
which a \Borel subset of a product of two \Polish spaces can be
separated from another subset by a countable union of rectangles
of low complexity was explored, as was the related question concerning
the existence of a minimal \Borel graph on a \Polish space with no
$\aleph_0$-coloring of low complexity. Positive results were obtained at
the second level of the \Borel hierarchy by adapting a classical proof of
\Hurewicz's original dichotomy theorem which, despite its simplicity,
nevertheless relies on a finite injury argument, and does not yield the
theorem in its natural generality. Moreover, the resulting minimal objects
in some sense fail to be canonical, in that they depend on an arbitrary
parameter. Positive results were also obtained at the third level of the
\Borel hierarchy, although the underlying arguments were quite intricate,
as were the corresponding minimal objects. The question as to whether
such results hold at the fourth level and beyond remains open. 

In \cite{JayneRogers}, it was shown that a function from an analytic
metric space to a separable metric space is \Gdelta-measurable if and
only if it is $\sigma$-continuous with closed witnesses. In \cite{Solecki},
this was shown to be a consequence of the fact that there is a
two-element basis, consisting of non-\Gdelta-measurable functions, for
the family of \Baire-class-one functions that are not
$\sigma$-continuous with closed witnesses. As a corollary, it was also
shown that if a \Baire-class-one function is not $\sigma$-continuous
with closed witnesses, then $\dominatingnumber$ is the least cardinal
$\kappa$ for which its domain is the union of $\kappa$-many closed
sets on which the function is continuous. These results were
established using ad-hoc recursive constructions reminiscent of those
behind the level-two \Lecomte-\Zeleny results. Despite having received
quite a bit of attention, the question as to whether these results
generalize in their most natural form to higher levels of the \Borel
hierarchy remains open.

Here we show that a common generalization of the
\Kechris-\Louveau-\Woodin generalization of \Hurewicz's original
dichotomy theorem (see \cite{KechrisLouveauWoodin}) and the
level-two results of \Lecomte-\Zeleny is a special case of the simplest
descriptive set-theoretic dichotomy theorem, as is the \Jayne-\Rogers
theorem. In addition to providing a simple unified explanation of
these phenomena and canonical minimal objects, this observation
also provides generalizations from analytic to separable metric spaces
under the axiom of determinacy, as well as the generalization of
\Solecki's result concerning $\dominatingnumber$ to \Borel functions.

In a future paper, we will use the particular structure of the minimal
objects so obtained to generalize \Solecki's basis theorem to all
\Borel functions, while simultaneously strengthening it using a finer
notion of embeddability. In another future paper, we will show that the
generalizations of \Hurewicz's theorem to every level of the difference
hierarchy at the second level of the \Borel hierarchy (including the
dichotomy theorem established in \cite{SaintRaymond} characterizing
the circumstances under which a \Borel set is a difference of two
\Fsigma sets) are also special cases of the simplest descriptive
set-theoretic dichotomy theorem. While our hope is that this approach
will lead to unified common generalizations to all levels of the \Borel
hierarchy, the question as to whether such arguments can even be
pushed to the third level remains open.

\subsection*{The underlying dichotomy}
The \definedterm{perfect set theorem} for a class $\Gamma$ of
topological spaces is the statement that for all $X \in \Gamma$, either
$\cardinality{X} \le \aleph_0$ or there is a continuous injection of
$\Cantorspace$ into $X$. In \cite{Souslin}, this simplest of descriptive
set-theoretic dichotomy theorems was established for analytic
\Hausdorff spaces. In \cite{Davis}, it was generalized to subsets of such
spaces under the axiom of determinacy.

The \definedterm{open graph dichotomy} for a class $\Gamma$ of
topological spaces is the statement that for all $X \in \Gamma$ and
open graphs $G$ on $X$, either $\chromaticnumber{G} \le \aleph_0$ or
there is a continuous homomorphism
from $\completegraph{\Cantorspace}$ to $G$. In \cite{Feng}, this
generalization of the perfect set theorem was established for analytic
\Hausdorff spaces, as was its generalization to subsets of such spaces
under the axiom of determinacy. This was achieved by showing that the
usual proofs of the perfect set theorem also yield the open graph
dichotomy, giving a sense in which the latter is also among the simplest
descriptive set-theoretic dichotomy theorems.

The \definedterm{box-open $D$-dimensional dihypergraph dichotomy}
for a class $\Gamma$ of topological spaces, or $\ogd[D]{\Gamma}$, is
the statement that for all $X \in \Gamma$ and box-open
$D$-dimensional dihypergraphs $H$ on $X$, either $\chromaticnumber
{H} \le \aleph_0$ or there is a continuous homomorphism from $\Heven
[D]$ to $H$. In \S\ref{ogd}, we establish this generalization of the open
graph dichotomy for analytic \Hausdorff spaces, as well as its
generalization to subsets of such spaces under the axiom of
determinacy. As in \cite{Feng}, we achieve this by showing that the
usual proofs of the perfect set theorem easily adapt to yield the open
dihypergraph dichotomy, giving a sense in which even the latter is
among the simplest descriptive set-theoretic dichotomy theorems as
well.

\subsection*{Applications}
In \S\ref{ksigma}, we give a first glimpse into how topological properties
can be codified into dihypergraphs by showing that the
\Kechris-\SaintRaymond generalizations of \Hurewicz's characterization
of the circumstances under which a \Polish space is \Ksigma (see \cite
{Kechris:Ksigma} or \cite{SaintRaymond:Ksigma}) is a special case of
the box-open dihypergraph dichotomy.

In \S\ref{compactifications}, we establish the basic properties of partial
compactifications of $\Bairespace$ arising from the box-open dihypergraph
dichotomy.

In \S\ref{hyperrectangles}, we show that a special case of the open
dihypergraph dichotomy yields a characterization of the circumstances
under which an analytic subset of a $D$-fold product of metric spaces
can be separated from another subset by a countable union of closed
hyperrectangles. As a corollary, we obtain a characterization of the
circumstances under which a $D$-dimensional dihypergraph on an
analytic metric space has a $\Deltaclass[0][2]$-measurable
$\aleph_0$-coloring. We also obtain the generalizations in which
analyticity is weakened to separability under the axiom of determinacy.
The \Kechris-\Louveau-\Woodin and level-two \Lecomte-\Zeleny
theorems are straightforward consequences of these results.

In \S\ref{sigmacontinuous}, we show that a special case of the box-open
dihypergraph dichotomy yields the \Jayne-\Rogers characterization
of the circumstances under which a function from an analytic metric
space to a separable metric space is $\sigma$-continuous with closed
witnesses. We also obtain the generalization in which analyticity is
weakened to separability under the axiom of determinacy.

In \S\ref{invariants}, we note that
$\chromaticnumber{\Heven[\N]} = \coveringnumber{\meagerideal}$ and
$\chromaticnumber{\Hevenprime[\N]} = \dominatingnumber$. In
conjunction with the box-open dihypergraph dichotomy, the latter fact easily
yields the promised generalization of \Solecki's theorem. We also
show that if $2 < D < \aleph_0$, then $\chromaticnumber
{\Heven[D]}$ is at least $\boundingnumber \cdot \coveringnumber
{\nullideal}$, consistently strictly below $\dominatingnumber$, and
consistently strictly above $\cofinality{\nullideal}$.

We work in the base theory $\ZF + \DC$ throughout, with the exception
of the final section, where we work in $\ZFC$ so as to keep our language
as transparent as possible.

\section{The box-open dihypergraph dichotomy} \label{ogd}

Despite the fact that it is essentially the same as the usual proof of the
perfect set theorem for analytic \Hausdorff spaces, we provide the proof
of the analogous instance of the open dihypergraph dichotomy for the
sake of the reader.

\begin{theorem} \label{ogd:dichotomy:analytic}
  Suppose that $D$ is a discrete space of cardinality at least two,
  $X$ is an analytic \Hausdorff space, and $H$ is a box-open
  $D$-dimensional dihypergraph on $X$. Then exactly one of the
  following holds:
  \begin{enumerate}
    \item There is an $\aleph_0$-coloring of $H$.
    \item There is a continuous homomorphism from $\Heven[D]$ to $H$.
  \end{enumerate}
\end{theorem}

\begin{theoremproof}
  To see that the two conditions are mutually exclusive, it is sufficient to
  show that there is no $\aleph_0$-coloring of $\Heven[D]$. Towards this
  end, suppose that $X \subseteq \functions{\N}{D}$ and $c \from X \to
  \N$ is an $\aleph_0$-coloring of $\restriction{\Heven[D]}{X}$,
  recursively find $d_n \in D$ such that $n \notin \image{c}{\extensions
  {\sequence{d_m}[m \le n]}}$ for all $n \in \N$, and observe that
  $\sequence{d_n}[n \in \N] \notin \preimage{c}{\N}$, thus $X \neq \functions{\N}{D}$.

  To see that at least one of the conditions hold, we can assume that
  $X \neq \emptyset$, in which case there is a continuous surjection
  $\pi \from \Bairespace \to X$. By replacing $H$ with its pullback
  through $\pi$, we can assume that $X = \Bairespace$.
  
  Set $S = \set{s \in \Bairetree}[\restriction{H}{\extensions{s}} \text{ has
  an $\aleph_0$-coloring}]$ and $Y = \setcomplement{\union[s \in S]
  [\extensions{s}]}$. Note that if $s \in \setcomplement{S}$, then there is
  no $\aleph_0$-coloring of $\restriction{H}{(\extensions{s} \intersection
  Y)}$, so there exists $\sequence{y_d}[d \in D] \in \restriction
  {H}{(\extensions{s} \intersection Y)}$, thus the fact that $H$ is box
  open yields a sequence $\sequence{s_d}[d \in D] \in \functions{D}
  {(\setcomplement{S})}$ of proper extensions of  $s$ such that
  $\product[d \in D][\extensions{s_d}] \subseteq H$. It follows that if
  there is no $\aleph_0$-coloring of $H$, then there is a function
  $f \from \functions{<\N}{D} \to \setcomplement{S}$ such that:
  \begin{enumerate}
    \renewcommand{\theenumi}{\alph{enumi}}
    \item $\forall d \in D \forall t \in \functions{<\N}{D} \ f(t)
      \strictlyextendedby f(t \concatenation \sequence{d})$.
    \item $\forall t \in \functions{<\N}{D} \ \product[d \in D][\extensions{f(t
      \concatenation \sequence{d})}] \subseteq H$.
  \end{enumerate}
  Condition (a) ensures that we obtain a continuous function $\phi \from
  \functions{\N}{D} \to Y$ by setting $\phi(\bfd) = \union[n \in \N]
  [f(\restriction{\bfd}{n})]$, and condition (b) implies that $\phi$ is a
  homomorphism from $\Heven[D]$ to $H$.
\end{theoremproof}

As in \cite{Feng}, the same argument yields the
natural generalization to $\kappa$-\Souslin \Hausdorff spaces, and a
derivative can be used to avoid the need for copious amounts of choice.
However, in order to establish the open dihypergraph dichotomy for all
subsets of analytic metric spaces under the axiom of determinacy, the
natural analog of the game considered in \cite[$\S{3}$]{Feng} is
insufficient, as it requires the first player to play $\cardinality{D}$-many
natural numbers in each round. We next show that a slowed-down
version of this game can be used instead.

Given a box-open $\N$-dimensional dihypergraph $H$ on
$\Bairespace$ and a set $X \subseteq \Bairespace$, consider the
$\omega$-length two-player game whose \textexponent{$n$}{th} round
consists of the first player playing a sequence $s_n \in \Bairetree$ with
the property that $i_m = 1 \implies s_m \strictlyextendedby s_n$ for all
$m < n$, and then the second playing a natural number $i_n < 2$.
There are two types of runs of the game, depending on whether the set
$N = \set{n \in \N}[i_n = 1]$ is finite or infinite. In the former case, the
first player wins if and only if $\product[n \in \N][\extensions{s_{m + n}}]
\subseteq H$, where $m$ is the least natural number strictly larger than
every element of $N$. In the latter, the first player wins if and only if
$\union[n \in N][s_n] \in X$.

\begin{proposition} \label{ogd:gameequivalent}
  Suppose that $H$ is a box-open $\N$-dimensional dihypergraph on
  $\Bairespace$ and $X \subseteq \Bairespace$.
  \begin{enumerate}
    \item The first player has a winning strategy if and only if there is a
      continuous homomorphism from $\Heven[\N]$ to $\restriction{H}
        {X}$.
    \item The second player has a winning strategy if and only if there is
      an $\aleph_0$-coloring of $\restriction{H}{X}$.
  \end{enumerate}
\end{proposition}
  
\begin{propositionproof}
  Let $S$ denote the set of sequences $s \in \Cantortree$ that do not
  end in zero. Let $\Ramseyspace$ denote the set of strictly increasing
  elements of $\Bairespace$, define $\pi_{\Bairespace, \Ramseyspace}
  \from \Bairespace \to \Ramseyspace$ by $\pi_{\Bairespace,
  \Ramseyspace}(b)(n) = n + \sum_{m \le n} b(m)$, and define $\pi_
  {\Ramseyspace, \Cantorspace} \from \Ramseyspace \to \Cantorspace$
  by $\pi_{\Ramseyspace, \Cantorspace}(b) = \characteristicfunction
  {\image{b}{\N}}$. Then the function $\pi_{\Bairespace, \Cantorspace} =
  \pi_{\Ramseyspace, \Cantorspace} \composition \pi_{\Bairespace,
  \Ramseyspace}$ is a homeomorphism from $\Bairespace$ to the
  space $C = \set{c \in \Cantorspace}[\cardinality{\support{c}} =
  \aleph_0]$, in addition to being an isomorphism of $\Heven[\N]$ with
  the restriction of the box-open $\N$-dimensional dihypergraph
  $H_{\Cantorspace} = \union[s \in S][{\product[n \in \N][\extensions{s
  \concatenation \constantsequence{0}{n} \concatenation \sequence
  {1}}]}]$ to $C$.
  
  Suppose that $\tau$ is a winning strategy for the first player, and
  define $\phi \from C \to \Bairespace$ by $\phi(c) = \union[n \in \support
  \negthickspace c][\tau(\restriction{c}{n})]$. To see that $\phi$ is a
  homomorphism from $\restriction{H_{\Cantorspace}}{C}$ to $H$, note
  that if $\sequence{c_n}[n \in \N] \in \restriction{H_{\Cantorspace}}{C}$,
  then there exists $s \in S$ for which $\sequence{c_n}[n \in \N] \in
  \product[n \in \N][\extensions{s \concatenation \constantsequence{0}{n}
  \concatenation \sequence{1}}]$, in which case $\sequence{\phi(c_n)}[n
  \in \N] \in \product[n \in \N][\extensions{\tau(s \concatenation
  \constantsequence{0}{n})}]$, so the fact that the first player wins runs
  of the game of the first type ensures that $\sequence{\phi(c_n)}[n \in
  \N] \in H$. The fact that the first player wins runs of the game of the
  second type implies that $\image{\phi}{C} \subseteq X$.
  
  Conversely, suppose that $\phi \from C \to X$ is a continuous
  homomorphism from $\restriction{H_{\Cantorspace}}{C}$ to $H$. For
  each non-empty sequence $s \in S$, let $\predecessor{s}$ denote the
  immediate predecessor of $s$. We will recursively construct functions
  $\sigma \from \Cantortree \to S$ and $\tau \from \Cantortree \to
  \Bairetree$ such that:
  \begin{enumerate}
    \renewcommand{\theenumi}{\alph{enumi}}
    \item $\forall s \in S \ \product[n \in \N][\extensions{\tau(s
      \concatenation \constantsequence{0}{n})}] \subseteq H$.
    \item $\forall s \in \Cantortree \ \image{\phi}{\extensions{\sigma(s)}}
      \subseteq \extensions{\tau(s)}$.
    \item $\forall n \in \N \forall s \in S \setminus \set{\emptysequence}
      \ \tau(\predecessor{s}) \strictlyextendedby \tau(s \concatenation
        \constantsequence{0}{n})$.
    \item $\forall n \in \N \forall s \in S \setminus \set{\emptysequence}
      \ \sigma(\predecessor{s})
        \concatenation \constantsequence{0}{n} \concatenation \sequence
          {1} \extendedby \sigma(s \concatenation \constantsequence{0}
            {n})$.
  \end{enumerate}
  Suppose that $s \in S$ and we have already found $\sigma(r)$ and
  $\tau(r)$ for all $r \strictlyextendedby s$. If $s = \emptysequence$,
  then set $u_s  = v_s = \emptysequence$, and otherwise define $u_s
  = \sigma(\predecessor{s})$ and $v_s = \tau(\predecessor{s})$. As
  $\sequence{u_s \concatenation \constantsequence{0}{n}
  \concatenation \constantsequence{1}{\infty}}[n \in \N] \in \restriction
  {H_{\Cantorspace}}{C}$, there are strict extensions $\tau(s
  \concatenation \constantsequence{0}{n}) \strictlyextendedby \phi(u_s
  \concatenation \constantsequence{0}{n} \concatenation
  \constantsequence{1}{\infty})$ of $v_s$ such that $\product[n \in \N]
  [\extensions{\tau(s \concatenation \constantsequence{0}{n})}]
  \subseteq H$, as well as positive integers $k_{n,s}$ such that $\image
  {\phi}{\extensions{u_s \concatenation \constantsequence{0}{n}
  \concatenation \constantsequence{1}{k_{n, s}}}} \subseteq \extensions
  {\tau(s \concatenation \constantsequence{0}{n})}$ for all $n \in \N$. We
  complete the construction by setting $\sigma(s \concatenation
  \constantsequence{0}{n}) = u_s \concatenation \constantsequence{0}
  {n} \concatenation \constantsequence{1}{k_{n,s}}$ for all $n \in \N$.
  To see that $\tau$ is a winning strategy for the first player, note that
  the first player wins runs of the game of the first type by condition (a),
  whereas the other conditions ensure that if $c \in C$, then $\phi
  (\union[n \in \support{c}][\sigma(\restriction{c}{n})]) = \union[n \in
  \support \negthickspace c][\tau(\restriction{c}{n})]$, thus the first player
  wins runs of the game of the second type.
  
  Suppose now that $\tau$ is a winning strategy for the second player.
  Let $T$ denote the set of partial runs of the game against $\tau$ for
  which $\sequence{\tau(\restriction{t}{\set{0, \ldots, n}})}[n < \length{t}]
  \in S$, and associate with each $t \in T$ the set $X_t = \set{x \in X}[t
  \neq \emptysequence \implies t(\length{t} - 1) \extendedby x]$. If
  $\sequence{x_n}[n \in \N] \in \restriction{H}{X_t}$, then there are
  sequences $t_n \strictlyextendedby x_n$ with the property that
  $\product[n \in \N][\extensions{t_n}] \subseteq H$ and $t \neq
  \emptysequence \implies \forall n \in \N \ t(\length{t} - 1)
  \strictlyextendedby t_n$, so the fact that the second player wins runs
  of the game of the first type therefore yields $n \in \N$ for which $t
  \concatenation \sequence{t_m}[m \le n] \in T$ and $x_n \in X_{t
  \concatenation \sequence{t_m}[m \le n]}$. In particular, it follows that
  the sets of the form $X_t \setminus \union[t \strictlyextendedby u, u \in
  T][X_u]$ are $H$-independent. As the fact that the second player wins
  runs of the game of the second type ensures that every $x \in X$
  appears in a set of the latter form, it follows that there is an
  $\aleph_0$-coloring of $\restriction{H}{X}$.
  
  Conversely, suppose that $c \from X \to \N$ is an $\aleph_0$-coloring
  of $\restriction{H}{X}$, and let $\tau$ be the strategy for the second
  player in which $0$ is played in the \textexponent{$n$}{th} round of the
  game if and only if $\preimage{c}{\set{k_n}} \intersection \extensions
  {s_n} \neq \emptyset$, where $k_n = \cardinality{\set{m < n}[i_m = 1]}$.
  The fact that the sets of the form $\preimage{c}{\set{k}}$ are
  $H$-independent ensures that the second player wins runs of the
  game of the first type, while the fact that $X \subseteq \preimage{c}
  {\N}$ implies that the second player wins runs of the game of the
  second type.
\end{propositionproof}

In particular, we obtain the following.

\begin{theorem}[$\AD$] \label{ogd:dichotomy:ad}
  Suppose that $Y$ is an analytic \Hausdorff space, $H$ is a box-open
  $\N$-dimensional dihypergraph on $Y$, and $X \subseteq Y$. Then
  exactly one of the following holds:
  \begin{enumerate}
    \item There is an $\aleph_0$-coloring of $\restriction{H}{X}$.
    \item There is a continuous homomorphism from $\Heven[\N]$ to
      $\restriction{H}{X}$.
  \end{enumerate}
\end{theorem}

\begin{theoremproof}
  As noted in the proof of Theorem \ref{ogd:dichotomy:analytic}, the two
  conditions are mutually exclusive, so it is sufficient to show that at
  least one of them holds. We can also assume that $Y \neq \emptyset$,
  in which case there is a continuous surjection $\phi \from \Bairespace
  \to Y$. By replacing $H$ and $X$ with their pullbacks through $\pi$, we
  can assume that $Y = \Bairespace$. But Proposition \ref
  {ogd:gameequivalent} then ensures that one of the two conditions holds.
\end{theoremproof}

\section{\Ksigma sets} \label{ksigma}

Given a topological space $X$, let $H_X$ denote the $\N$-dimensional
dihypergraph on $X$ consisting of all injective sequences $\sequence
{x_n}[n \in \N]$ of elements of $X$ with no convergent subsequence.

\begin{proposition} \label{ksigma:hypergraph}
  Suppose that $X$ is a metric space.
  \begin{enumerate}
    \item The dihypergraph $H_X$ is box open.
    \item There is an $\aleph_0$-coloring of the restriction of $H_X$ to a
      set $Y \subseteq X$ if and only if $Y$ is contained in a \Ksigma
      subset of $X$.
    \item A continuous function $\phi \from \Bairespace \to X$ is a
      homomorphism from $\Heven[\N]$ to $H_X$ if and only if it is an
      injective closed map.
  \end{enumerate}
\end{proposition}

\begin{propositionproof}
  To see (1), note that if $\sequence{x_n}[n \in \N] \in H_X$, then there
  exist positive real numbers $\epsilon_n \goesto 0$ such that $\rho_X
  (x_m, x_n) \ge 2\epsilon_n$ for all natural numbers $m \neq n$, in
  which case $\product[n \in \N][{\ball[X]{x_n}{\epsilon_n}}] \subseteq H_X$.
  
  To see (2), it is sufficient to observe that a set $Y \subseteq X$ is
  $H_X$-independent if and only if its closure is compact.
  
  To see (3), note first that if $\phi$ is an injective closed map, then the
  fact that each sequence $\sequence{\bfd_n}[n \in \N] \in \Heven[\N]$
  is an injective enumeration of a closed discrete set ensures that the
  same holds of $\sequence{\phi(\bfd_n)}[n \in \N]$. Conversely,
  suppose that $\phi$ is a homomorphism from $\Heven[\N]$ to $H_X$.
  The fact that $H_X$ consists solely of injective sequences easily
  implies that $\phi$ is injective. To see that $\phi$ is a closed map, it is
  sufficient to show that every sequence $\sequence{\bfd_n}[n \in \N]$
  of elements of $\Bairespace$ for which $\sequence{\phi(\bfd_n)}[n \in
  \N]$ converges has a convergent subsequence. If there exists $\bfd
  \in \Bairespace$ such that $\bfd_n(i) < \bfd(i)$ for all
  $i, n \in \N$, then the compactness of $\product[i \in \N][\bfd(i)]$ yields
  the desired subsequence. So suppose, towards a contradiction, that
  there does not exist such a $\bfd$. Then there is a least $k
  \in \N$ for which $\set{\bfd_n(k)}[n \in \N]$ is infinite. By passing to a
  subsequence, we can assume that for all distinct $m, n \in \N$, the
  sequences $\bfd_m$ and $\bfd_n$ differ from one another for the first
  time on their \textexponent{$k$}{th} coordinates. By passing to a further
  subsequence, we can assume that $\sequence{\bfd_n}[n \in \N]$ is a
  subsequence of an element of $\Heven[\N]$, so $\sequence{\phi
  (\bfd_n)}[n \in \N]$ is a subsequence of an element of $H_X$,
  contradicting the fact that it converges.
\end{propositionproof}

In particular, we obtain the following.

\begin{theorem}[{$\ogd[\N]{\Gamma}$}]
  \label{ksigma:dichotomy}
  Suppose that $X$ is a metric space and $Y \subseteq X$ is in
  $\Gamma$. Then exactly one of the following holds:
  \begin{enumerate}
    \item The set $Y$ is contained in a \Ksigma subset of $X$.
    \item There is a closed continuous injection $\phi \from
      \Bairespace \to X$ with the property that $\image{\phi}
        {\Bairespace} \subseteq Y$.
  \end{enumerate}
\end{theorem}

\begin{theoremproof}
  This follows from Proposition \ref{ksigma:hypergraph}.
\end{theoremproof}

The special cases of Theorem \ref{ksigma:dichotomy} where
$\Gamma$ is either the family of analytic spaces or the family of
separable spaces easily yields the \Kechris-\SaintRaymond
generalizations of \Hurewicz's characterization of the circumstances
under which a \Polish space is \Ksigma.

\section{Partial compactifications} \label{compactifications}

For each topological space $X$ and discrete set $D \subseteq
X$, endow the set $\modifiedpower{\N}{D}{X} = \functions{\N}{D} \union
\set{t \concatenation \sequence{x}}[t \in \functions{<\N}{D} \mathand x
\in \protect\setcomplement{D}]$ with the topology generated by the
sets $\extensions{t}[U] = \set{\bfd \in \modifiedpower{\N}{D}{X}}[t
\strictlyextendedby \bfd \mathand \bfd(\length{t}) \in U]$, where $t \in
\functions{<\N}{D}$ and $U \subseteq X$ is open. Such spaces arise
naturally in applications of the open dihypergraph dichotomy. An
instructive example to consider is $\modifiedpower{\N}{\N}
{\compactification{\N}}$, where $\compactification{\N} = \N \union \set
{\infty}$ denotes the one-point compactification of $\N$.

\begin{proposition} \label{compactification:compact}
  Suppose that $X$ is a compact space and $D \subseteq X$ is
  discrete and open. Then $\modifiedpower{\N}{D}{X}$ is compact.
\end{proposition}

\begin{propositionproof}
  Set $\calU = \set{\extensions{t}[U]}[t \in \functions{<\N}{D} \mathand U
  \subseteq X \text{ is open}]$, and for each sequence $t \in \functions
  {<\N}{D}$ and family $\calV \subseteq \calU$, let $\calV_t$ denote the
  family of open sets $V \subseteq X$ for which $\extensions{t}[V] \in
  \calV$.
  
  \begin{lemma} \label{compactification:compact:extension}
    Suppose that $t \in \functions{<\N}{D}$ and $\calV \subseteq \calU$
    covers $\extensions{t}[X]$. If for all $d \in D$ there is a finite set
    $\calF_d \subseteq \calV$ covering $\extensions{t \concatenation
    \sequence{d}}[X]$, then there is a finite set $\calF \subseteq \calV$
    covering $\extensions{t}[X]$.
  \end{lemma}
  
  \begin{lemmaproof}
    We can assume that there do not exist $s \strictlyextendedby
    t$ and an open set $U \subseteq X$ for which $t(\length{s}) \in U$ and
    $\extensions{s}[U] \in \calV$. Then there is a finite set $\calF
    \subseteq \calV_t$ covering $\setcomplement{D}$, in which case the
    set $F = \setcomplement{\union[\calF]}$ is compact and contained in
    $D$, thus finite. But $\calF \union \union[d \in F][\calF_d]$
    covers $\extensions{t}[X]$.
  \end{lemmaproof}
  
  If $\calV \subseteq \calU$ is a cover of $\modifiedpower{\N}{D}{X}$
  with no finite subcover, then by recursively applying the contrapositive
  of Lemma \ref{compactification:compact:extension}, we obtain a
  sequence $\bfd \in \functions{\N}{D}$ such that for no $n \in \N$ is
  there a finite subset of $\calV$ covering $\extensions{\restriction{\bfd}
  {n}}[X]$, contradicting the fact that $\bfd \in \union[\calV]$.
\end{propositionproof}

\begin{proposition}
  Suppose that $X$ is a (complete) ultrametric space and $D \subseteq
  X$ is discrete and open. Then $\modifiedpower{\N}{D}{X}$ has a
  (complete) compatible ultrametric.
\end{proposition}

\begin{propositionproof}
  By replacing $\rho_X$ with $\rho_X / (1 + \rho_X)$, we can assume
  that $\rho_X < 1$. Fix real numbers $\epsilon_d > 0$ such that
  $\sup_{d \in D} \epsilon_d < 1$ and $\ball{d}{\epsilon_d} = \set{d}$ for
  all $d \in D$, and define $\rho \from \modifiedpower{\N}{D}{X} \times
  \modifiedpower{\N}{D}{X} \to \R$ by
  \begin{equation*}
    \rho(\bfc, \bfd) =
        \begin{cases}
          0 & \text{if $\bfc = \bfd$ and} \\
          \rho_X(\bfc(n), \bfd(n)) \product[m < n][\epsilon_{\bfc(m)}] &
            \text{otherwise,}
        \end{cases}
  \end{equation*}
  where $n = n(\bfc, \bfd)$ is the least natural number for which $\bfc(n)
  \neq \bfd(n)$.
  
  To establish that $\rho$ is an ultrametric, it is sufficient to show that if
  $\bfb, \bfc, \bfd \in \modifiedpower{\N}{D}{X}$ are distinct, then $\rho
  (\bfb, \bfd) \le \max \set{\rho(\bfb, \bfc), \rho(\bfc, \bfd)}$. Setting $n
  = \max \set{n(\bfb, \bfc), n(\bfc, \bfd)}$, there are three cases to check:
  \begin{enumerate}
    \item If $n(\bfb, \bfd) < n$, then $\rho(\bfb, \bfd) \in \set{\rho(\bfb, \bfc),
      \rho(\bfc, \bfd)}$.
    \item If $n(\bfb, \bfd) = n$, then $n(\bfb, \bfc) = n(\bfb, \bfd) = n(\bfc,
    \bfd)$, so
      \begin{align*}
        \rho(\bfb, \bfd)
          & \textstyle = \rho_X(\bfb(n), \bfd(n)) \product[m < n]
            [\epsilon_{\bfb(m)}] \\
          & \textstyle \le \max \set{\rho_X(\bfb(n), \bfc(n)), \rho_X(\bfc(n),
            \bfd(n))} \product[m < n][\epsilon_{\bfb(m)}] \\
          & = \max \set{\rho(\bfb, \bfd), \rho(\bfc, \bfd)}.
      \end{align*}
    \item If $n(\bfb, \bfd) > n$, then $\rho(\bfb, \bfd) < \product[m \le n]
      [\epsilon_{\bfb(m)}] \le \rho(\bfb, \bfc)$.
  \end{enumerate}

  To see that the topology generated by $\rho$ is contained in that
  of $\modifiedpower{\N}{D}{X}$, suppose that $\bfd \in \modifiedpower
  {\N}{D}{X}$ and $\epsilon > 0$. If $\bfd \in \functions{\N}{D}$, then
  there exists $n \in \N$ sufficiently large that $\product[m < n][\epsilon_
  {\bfd(m)}] \le \epsilon$, in which case $\extensions{\restriction{\bfd}{n}}
  [X]$ is an open neighborhood of $\bfd$ contained in $\ball{\bfd}
  {\epsilon}$. Otherwise, there exist $t \in \functions{<\N}{D}$ and $x \in
  \setcomplement{D}$ for which $\bfd = t \concatenation \sequence{x}$,
  in which case $\extensions{t}[{\ball[X]{x}{\epsilon}}]$ is an open
  neighborhood of $\bfd$ contained in $\ball{\bfd}{\epsilon}$.

  To see that the topology of $\modifiedpower{\N}{D}{X}$ is contained in
  that generated by $\rho$, suppose that $t \in \functions{<\N}{D}$, $U
  \subseteq X$ is open, and $\bfd \in \extensions{t}[U]$, and fix $0 <
  \epsilon < 1$ such that $\ball[X]{\bfd(\length{t})}{\epsilon} \subseteq U$.
  Then $\ball{\bfd}{\epsilon \product[n < \length{t}][\epsilon_{t(n)}]}
  \subseteq \extensions{t}[U]$.
  
  To see that the completeness of $\rho_X$ yields that of $\rho$,
  suppose that $\sequence{\bfd_k}[k \in \N]$ is an injective \Cauchy
  sequence of elements of $\modifiedpower{\N}{D}{X}$. If there exists
  $\bfd \in \functions{\N}{D}$ with the property that $\bfd_k(n) \goesto
  \bfd(n)$ for all $n \in \N$, then $\bfd_k \goesto \bfd$. Otherwise, there
  is a sequence $t \in \functions{<\N}{D}$ of maximal length such that
  $\bfd_k(n) \goesto t(n)$ for all $n < \length{t}$. By passing to a
  subsequence of $\sequence{\bfd_k}[k \in \N]$, we can assume that
  $n(\bfd_j, \bfd_k) = \length{t}$, thus $\rho(\bfd_j, \bfd_k) = \rho_X
  (\bfd_j(\length{t}), \bfd_k(\length{t})) \product[n < \length{t}]
  [\epsilon_{t(n)}]$, for all distinct $j, k \in \N$. It follows that $\sequence
  {\bfd_k(\length{t})}[k \in \N]$ is \Cauchy, and therefore convergent,
  thus $\sequence{\bfd_k}[k \in \N]$ is convergent as well.
\end{propositionproof}

Let $\Rcon{X}$ denote the set of convergent sequences $\sequence
{x_n}[n \in \N]$ of elements of $X$, and let $\Rcon[D]{X}$ denote the
$\N$-ary relation on $\functions{\N}{D}$ given by $\Rcon[D]{X} = \union
[\pair{\bfd}{t} \in (\Rcon{X} \intersection \functions{\N}{D}) \times
\functions{<\N}{D}][{\product[n \in \N][\extensions{t \concatenation
\sequence{\bfd(n)}}]}]$. In our applications of the open dihypergraph
dichotomy, the following fact will yield extensions of homomorphisms
between dihypergraphs.

\begin{proposition} \label{compactification:extension}
  Suppose that $X$ and $Y$ are metric spaces, $D \subseteq X$ is
  dense, discrete, and open, and $\phi \from \functions{\N}{D} \to Y$ is
  a continuous homomorphism from $\Rcon[D]{X}$ to $\Rcon{Y}$.
  Then there is a continuous extension of $\phi$ to $\modifiedpower{\N}
  {D}{X}$.
\end{proposition}

\begin{propositionproof}
  Given a point $x \in X$, we say that a sequence $\sequence{X_n}[n
  \in \N]$ of subsets of $X$ \definedterm{converges} to $x$, or $X_n
  \goesto x$, if for every open neighborhood $U \subseteq X$ of $x$,
  there exists $n \in \N$ for which $\union[m \ge n][X_m] \subseteq U$.

  \begin{lemma} \label{compactification:extension:convergence}
    Suppose that $t \in \functions{<\N}{D}$ and $x \in \setcomplement
    {D}$. Then there exists $y_{t,x} \in Y$ such that $\bfd(n) \goesto x
    \implies \image{\phi}{\extensions{t \concatenation \sequence{\bfd
    (n)}}} \goesto y_{t,x}$ for all $\bfd \in \functions{\N}{D}$.
  \end{lemma}
  
  \begin{lemmaproof}
    Note that if $\bfd \in \functions{\N}{D}$, $\bfd(n) \goesto x$, and
    $y_n \in \image{\phi}{\extensions{t \concatenation \sequence{\bfd
    (n)}}}$ for all $n \in \N$, then there exists $y_{t,x} \in Y$ for which
    $y_n \goesto y_{t,x}$. If there exists $\bfc \in \functions{\N}{D}$ such
    that $\bfc(n) \goesto x$ and $\image{\phi}{\extensions{t
    \concatenation \sequence{\bfc(n)}}} \doesnotgoto y_{t,x}$, then there
    exist an infinite set $N \subseteq \N$, an open neighborhood $V
    \subseteq Y$ of $y_{t,x}$, and points $y_n \in \image{\phi}
    {\extensions{t \concatenation \sequence{\bfc(n)}}} \setminus V$ for
    all $n \in N$. By thinning down $N$, we can assume that it is
    co-infinite. Set $\bfb = (\restriction{\bfc}{N}) \union (\restriction{\bfd}
    {\setcomplement{N}})$, and observe that $\bfb(n) \goesto x$ but
    $\sequence{\phi(\bfb(n))}[n \in \N]$ does not converge, a
    contradiction.
  \end{lemmaproof}
  
  To see that the extension given by $\extension{\phi}(t \concatenation
  \sequence{x}) = y_{t,x}$ is continuous, suppose that $\bfd \in
  \modifiedpower{\N}{D}{X}$ and $V \subseteq Y$ is an open
  neighborhood of $\extension{\phi}(\bfd)$, and fix an open
  neighborhood $W \subseteq Y$ of $\extension{\phi}(\bfd)$ whose
  closure is contained in $V$. If $\bfd \in \functions{\N}{D}$, then the
  continuity of $\phi$ yields $n \in \N$ for which $\image{\phi}
  {\extensions{\restriction{\bfd}{n}}} \subseteq W$, in which case
  $\image{\extension{\phi}}{\extensions{\restriction{\bfd}{n}}[X]} \subseteq
  \closure{\image{\phi}{\extensions{\restriction{\bfd}{n}}}} \subseteq
  \closure{W} \subseteq V$. Otherwise, there exists $t \in \functions{<\N}
  {D}$ and $x \in \setcomplement{D}$ for which $\bfd = t \concatenation
  \sequence{x}$, so Lemma \ref
  {compactification:extension:convergence} yields an open
  neighborhood $U \subseteq X$ of $x$ for which $\image{\phi}
  {\extensions{t}[U] \intersection \functions{\N}{D}} \subseteq W$, in
  which case $\image{\extension{\phi}}{\extensions{t}[U]} \subseteq
  \closure{\image{\phi}{\extensions{t}[U]}} \subseteq \closure{W}
  \subseteq V$.
\end{propositionproof}

\section{Countable unions of hyperrectangles} \label{hyperrectangles}

A \definedterm{hyperrectangular homomorphism} from a sequence
$\sequence{R_i}[i \in I]$ of subsets of $\product[d \in D][X_d]$ to a
sequence $\sequence{S_i}[i \in I]$ of subsets of $\product[d \in D]
[Y_d]$ is a function $\phi$ of the form $\product[d \in D][\phi_d]$,
where $\phi_d \from \union[i \in I][\image{\projection[d]}{R_i}] \to Y_d$
for all $d \in D$, such that $\image{\phi}{R_i} \subseteq S_i$ for all $i
\in I$.

We use $\Hcountable[D]$ to denote the $D$-dimensional dihypergraph
on $\modifiedpower{\N}{(D \times \N)}{D \times \compactification{\N}}$
consisting of all sequences $\sequence{t \concatenation \sequence
{\pair{d}{\infty}}}[d \in D]$, where $t \in \functions{<\N}{(D \times \N)}$.

Given a sequence $\sequence{X_d}[d \in D]$ and disjoint sets $R, S
\subseteq \product[d \in D][X_d]$, let $H_{R,S}$ denote the $(D \times
\N)$-dimensional dihypergraph on $R$ consisting of all sequences
$\sequence{\sequence{x_{c,d,n}}[c \in D]}[\pair{d}{n} \in D \times \N]$
of elements of $R$ with the property that $\overline{x}_d = \lim_{n
\goesto \infty} x_{d,d,n}$ exists for all $d \in D$ and $\sequence
{\overline{x}_d}[d \in D] \in S$.
  
\begin{proposition} \label{hyperrectangles:hypergraph}
  Suppose that $D$ is a discrete set, $\sequence{X_d}[d \in D]$ is a 
  sequence of metric spaces, and $R, S \subseteq \product[d \in D]
  [X_d]$ are disjoint.
  \begin{enumerate}
    \item The $(D \times \N)$-dimensional dihypergraph $H_{R, S}$ is
      box open.
    \item There is an $\aleph_0$-coloring of $H_{R, S}$ if and only if
      there is a countable union of closed hyperrectangles separating
      $R$ from $S$.
    \item There is a continuous homomorphism from $\Heven[(D \times
      \N)]$ to $H_{R, S}$ if and only if there is a continuous
      hyperrectangular homomorphism from $\pair{\diagonal[D]
      {\functions{\N}{(D \times \N)}}}{\Hcountable[D]}$ to $\pair{R}{S}$.
  \end{enumerate}
\end{proposition}

\begin{propositionproof}
  To see (1), note that if $\sequence{\sequence{x_{c,d,n}}[c \in D]}[\pair
  {d}{n} \in D \times \N] \in H_{R, S}$, $\epsilon_n \goesto 0$, and
  $U_{d,n} = \set{\sequence{x_c}[c \in D] \in R}[\rho_{X_d}(x_d, x_{d,d,n}) <
  \epsilon_n]$ for all $\pair{d}{n} \in D \times \N$, then $\product[\pair{d}
  {n} \in D \times \N][U_{d,n}] \subseteq H_{R, S}$.
  
  To see (2), note first that if $Q \subseteq R$ and $\sequence{\overline
  {x}_d}[d \in D] \in \product[d \in D][\closure{\image{\projection[d]}{Q}}]$,
  then there are sequences $\sequence{x_{d,n}}[n \in \N]$ of elements
  of $\image{\projection[d]}{Q}$ such that $x_{d,n} \goesto \overline
  {x}_d$ for all $d \in D$, so there are sequences $\sequence{x_{c,d,n}}
  [c \in D] \in Q$ such that $x_{d,d,n} = x_{d,n}$ for all $\pair{d}{n} \in D
  \times \N$, thus $x_{d,d,n} \goesto \overline{x}_d$ for all $d \in D$. It
  follows that if $Q$ is $H_{R,S}$-independent, then $\product[d \in D]
  [\closure{\image{\projection[d]}{Q}}]$ and $S$ are disjoint, so if $c
  \from X \to \N$ is an $\aleph_0$-coloring of $H_{R,S}$, then the union
  of the closed hyperrectangles $\product[d \in D][\closure{\image
  {\projection[d]}{\preimage{c}{\set{n}}}}]$ separates $R$ from $S$.
  Conversely, suppose that $F_d$ is a closed subset of $X_d$ for all $d
  \in D$ and $S \intersection \product[d \in D][F_d] = \emptyset$. If
  $\sequence{\sequence{x_{c,d,n}}[c \in D]}[\pair{d}{n} \in D \times \N]$
  is a sequence of elements of $\product[c \in D][F_c]$ such that
  $\overline{x}_d = \lim_{n \goesto \infty} x_{d,d,n}$ exists for all $d \in
  D$, then $\sequence{\overline{x}_d}[d \in D] \in \product[d \in D][F_d]$,
  so $\sequence{\overline{x}_d}[d \in D] \notin S$, thus $R \intersection
  \product[d \in D][F_d]$ is $H_{R, S}$-independent. Hence if there is a
  countable union of closed hyperrectangles separating $R$ from $S$,
  then there is an $\aleph_0$-coloring of $H_{R,S}$.

  To see (3), suppose first that $\product[d \in D][\phi_d]$ is a
  continuous hyperrectangular homomorphism from $\pair{\diagonal[D]
  {\functions{\N}{(D \times \N)}}}{\Hcountable[D]}$ to $\pair{R}{S}$, and
  define $\phi \from \functions{\N}{(D \times \N)} \to \product[d \in D]
  [X_d]$ by $\phi(\bfd) = \sequence{\phi_d(\bfd)}[d \in D]$. Clearly $\phi$
  is continuous and $\image{\phi}{\functions{\N}{(D \times \N)}} = \image
  {(\product[d \in D][\phi_d])}{\diagonal[D]{\functions{\N}{(D \times \N)}}}
  \subseteq R$. To see that $\phi$ is a homomorphism from $\Heven[(D
  \times \N)]$ to $H_{R, S}$, we need only note that $\sequence{\lim_{n
  \goesto \infty} \phi_d(\extensions{t \concatenation \sequence{\pair{d}
  {n}}})}[d \in D] = \sequence{\phi_d(t \concatenation \sequence{\pair{d}
  {\infty}})}[d \in D]$ for all $t \in \functions{<\N}{(D \times \N)}$, since
  $\image{(\product[d \in D][\phi_d])}{\Hcountable[D]} \subseteq S$.
  Conversely, suppose that $\phi \from \functions{\N}{(D \times \N)} \to
  R$ is a continuous homomorphism from $\Heven[(D \times \N)]$ to
  $H_{R,S}$. Then each of the functions $\projection[d] \composition
  \phi$ is a continuous homomorphism from $\Rcon[(D \times \N)]{(D
  \times \N) \union \set{\pair{d}{\infty}}}$ to $\Rcon{X_d}$, so Proposition
  \ref{compactification:extension} yields continuous extensions $\phi_d
  \from \modifiedpower{\N}{(D \times \N)}{(D \times \N) \union \set{\pair
  {d}{\infty}}} \to X_d$ for all $d \in D$, and clearly $\product[d \in D]
  [\phi_d]$ is a hyperrectangular homomorphism from $\pair{\diagonal
  [D]{\functions{\N}{(D \times \N)}}}{\Hcountable[D]}$ to $\pair{R}{S}$.
\end{propositionproof}

In particular, we obtain the following.

\begin{theorem}[{$\ogd[\N]{\Gamma}$}]
  \label{hyperrectangles:dichotomy}
  Suppose that $D$ is a countable discrete set, $\sequence{X_d}[d \in
  D]$ is a sequence of metric spaces, $R, S \subseteq \product[d \in D]
  [X_d]$ are disjoint, and $R \in \Gamma$. Then exactly one of the
  following holds:
  \begin{enumerate}
    \item There is a countable union of closed hyperrectangles
      separating $R$ from $S$.
    \item There is a continuous hyperrectangular
      homomorphism from $\pair{\diagonal[D]{\functions{\N}{(D \times
      \N)}}}{\Hcountable[D]}$ to $\pair{R}{S}$.
  \end{enumerate}
\end{theorem}

\begin{theoremproof}
  This follows from Proposition \ref{hyperrectangles:hypergraph}.
\end{theoremproof}

The special cases of Theorem \ref{hyperrectangles:dichotomy} where
$D = 1$ and $\Gamma$ is either the family of analytic spaces or the
family of separable spaces easily yield the \Kechris-\Louveau-\Woodin
generalizations of \Hurewicz's characterization of the circumstances
under which two sets can be separated by an \Fsigma set.

The special case of Theorem \ref{hyperrectangles:dichotomy}
where $D = 2$ and $\Gamma$ is the family of analytic spaces
easily yields the \Lecomte-\Zeleny characterization of the circumstances
under which an analytic subset of the plane can be separated from another
subset of the plane by a countable union of closed rectangles. The special
case where $D = 2$ and $\Gamma$ is the family of separable spaces
easily yields the generalization in which analyticity is weakened to
separability under the axiom of determinacy.

\begin{proposition} \label{hyperrectangles:hypergraphs}
  Suppose that $D$ is a discrete set of cardinality at least two, $X$ is a
  metric space, and $H$ is a $D$-dimensional dihypergraph on $X$.
  \begin{enumerate}
    \item There is a countable union of closed hyperrectangles separating
      $\diagonal[D]{X}$ from $H$ if and only if there is a $\Deltaclass[0]
      [2]$-measurable $\aleph_0$-coloring of $H$.
    \item There is a continuous hyperrectangular homomorphism from
      $\pair{\diagonal[D]{\functions{\N}{(D \times \N)}}}{\Hcountable[D]}$
      to $\pair{\diagonal[D]{X}}{H}$ if and only if there is a continuous
      homomorphism from $\Hcountable[D]$ to $H$.
  \end{enumerate}
\end{proposition}

\begin{propositionproof}
  To see (1), note first that if $\sequence{F_{d,n}}[\pair{d}{n} \in D \times
  \N]$ is a sequence of closed subsets of $X$ for which $\union[n \in \N]
  [{\product[d \in D][F_{d,n}]}]$ separates $\diagonal[D]{X}$ from $H$,
  then the $H$-independent closed sets $\intersection[d \in D][F_{d,n}]$
  cover $X$, so there is a $\Deltaclass[0][2]$-measurable
  $\aleph_0$-coloring of $H$. Conversely, if there is a $\Deltaclass[0]
  [2]$-measurable $\aleph_0$-coloring of $H$, then there is a cover
  $\sequence{F_n}[n \in \N]$ of $X$ by $H$-independent closed sets,
  in which case $\union[n \in \N][\functions{D}{F_n}]$ is a
  countable union of closed hyperrectangles separating $\diagonal[D]
  {X}$ from $H$.
  
  To see (2), note first that if $\phi \from \modifiedpower{\N}{(D \times
  \N)}{D \times \compactification{\N}} \to X$ is a homomorphism from
  $\Hcountable[D]$ to $H$, and $\phi_d$ is the restriction of $\phi$ to
  $\modifiedpower{\N}{(D \times \N)}{(D \times \N) \union \set{\pair
  {d}{\infty}}}$ for all $d \in D$, then $\product[d \in D][\phi_d]$ is a
  hyperrectangular homomorphism from $\pair{\diagonal[D]{\functions
  {\N}{(D \times \N)}}}{\Hcountable[D]}$ to $\pair{\diagonal[D]{X}}{H}$.
  Conversely, observe that if $\product[d \in D][\phi_d]$ is a continuous
  hyperrectangular homomorphism from $\pair{\diagonal[D]{\functions
  {\N}{(D \times \N)}}}{\Hcountable[D]}$ to $\pair{\diagonal[D]{X}}{H}$,
  then the function $\phi = \union[d \in D][\phi_d]$ is a continuous
  homomorphism from $\Hcountable[D]$ to $H$.
\end{propositionproof}

In particular, we obtain the following.

\begin{theorem}[{$\ogd[\N]{\Gamma}$}]
  \label{hyperrectangles:dichotomy:hypergraphs}
  Suppose that $D$ is a countable discrete set of cardinality at least two,
  $X$ is a metric space in $\Gamma$, and $H$ is a $D$-dimensional
  dihypergraph on $X$. Then exactly one of the following holds:
  \begin{enumerate}
    \item There is a $\Deltaclass[0][2]$-measurable $\aleph_0$-coloring
      of $H$.
    \item There is a continuous homomorphism from $\Hcountable[D]$ to
      $H$.
  \end{enumerate}
\end{theorem}

\begin{theoremproof}
  This follows from Propositions \ref{hyperrectangles:hypergraph} and
  \ref{hyperrectangles:hypergraphs}.
\end{theoremproof}

The special case of Theorem \ref
{hyperrectangles:dichotomy:hypergraphs} where $D = 2$ and
$\Gamma$ is the family of analytic spaces easily yields the
\Lecomte-\Zeleny characterization of the circumstances under which
a graph on an analytic metric space has a $\Deltaclass[0][2]$-measurable
$\aleph_0$-coloring. The special case where $D = 2$ and $\Gamma$ is
the family of separable spaces easily yields the generalization in which
analyticity is weakened to separability under the axiom of determinacy.

\section{Sigma-continuous functions with closed witnesses} \label
  {sigmacontinuous}

A \definedterm{reduction} of a set $R \subseteq \functions{D}{X}$ to a
set $S \subseteq \functions{D}{Y}$ is a single function $\phi \from X \to
Y$ that is a homomorphism from $R$ to $S$ and from $\setcomplement
{R}$ to $\setcomplement{S}$.

Given a function $\pi \from X \to Y$, let $H_\pi$ denote the
$\N$-dimensional dihypergraph on $\graph{\pi}$ consisting of all
sequences $\sequence{\pair{x_n}{y_n}}[n \in \N]$ of elements of
$\graph{\pi}$ with the property that $\overline{x} = \lim_{n \goesto \infty}
x_n$ exists but $\pi(\overline{x}) \notin \closure{\set{y_n}[n \in \N]}$.

\begin{proposition}
  \label{sigmacontinuous:hypergraph}
  Suppose that $X$ and $Y$ are metric spaces and $\pi \from X \to
  \nolinebreak Y$.
  \begin{enumerate}
    \item The dihypergraph $H_\pi$ is box open.
    \item There is an $\aleph_0$-coloring of $H_\pi$ if and only if $\pi$
      is $\sigma$-continuous with closed witnesses.
    \item There is a continuous homomorphism from $\Heven[\N]$ to
      $H_\pi$ if and only if there is a continuous function $\psi \from
      \modifiedpower{\N}{\N}{\compactification{\N}} \to X$ with the
      property that $\pi \composition \psi$ is a reduction of $\Bairespace$
      to a closed set and $\restriction{(\pi \composition \psi)}
      {\Bairespace}$ is continuous.
  \end{enumerate}
\end{proposition}

\begin{propositionproof}
  To see (1), note that if $\sequence{\pair{x_n}{y_n}}[n \in \N] \in H_\pi$
  and $\overline{x} = \lim_{n \goesto \infty} x_n$, then there exists
  $\epsilon > 0$ such that $\rho_Y(\pi(\overline{x}), y_n) \ge 2\epsilon$ for
  all $n \in \N$. Fix positive real numbers $\epsilon_n \goesto 0$, and
  observe that the intersection of $\graph{\pi}$ with $\product[n \in \N]
  [{\ball[X]{x_n}{\epsilon_n} \times \ball[Y]{y_n}{\epsilon}}]$ is contained
  in $H_\pi$.
  
  To see (2), note first that if $F \subseteq X$ is a closed set on which
  $\pi$ is continuous, then $\graph{\pi} \intersection (F \times Y)$ is
  $H_\pi$-independent, so if $\pi$ is $\sigma$-continuous with closed
  witnesses, then there is a \Borel $\aleph_0$-coloring of $H_\pi$.
  Conversely, if $R \subseteq \graph{\pi}$ and $\sequence
  {\overline{x}_n}[n \in \N]$ is a convergent sequence of elements of the
  closure of $\image{\projection[X]}{R}$, then there are sequences
  $\sequence{x_{m,n}}[m \in \N]$ of elements of $\image{\projection[X]}
  {R}$ such that $x_{m,n} \goesto \overline{x}_n$ for all $n \in \N$. If
  $R$ is $H_\pi$-independent, then $\pi(x_{m,n}) \goesto \pi(\overline
  {x}_n)$ for all $n \in \N$, so there is a function $f \from \N \to \N$ such
  that $\rho_X(x_{f(n), n}, \overline{x}_n) \goesto 0$ and $\rho_Y(\pi
  (x_{f(n), n}), \pi(\overline{x}_n)) \goesto 0$. It follows that $\sequence
  {x_{f(n), n}}[n \in \N]$ converges, so the fact that $R$ is
  $H_\pi$-independent ensures that $\pi(\lim_{n \goesto \infty} x_{f(n),
  n}) = \lim_{n \goesto \infty} \pi(x_{f(n), n})$, in which case
  $\pi(\lim_{n \goesto \infty} \overline{x}_n) = \lim_{n \goesto \infty} \pi
  (\overline{x}_n)$, thus $\restriction{\pi}{\closure{\image{\projection[X]}
  {R}}}$ is continuous. In particular, it follows that if there is an
  $\aleph_0$-coloring of $H_\pi$, then $\pi$ is $\sigma$-continuous
  with closed witnesses.
  
  To see (3), note first that if $\psi \from \modifiedpower{\N}{\N}
  {\compactification{\N}} \to X$ is continuous, $\pi \composition \psi$ is
  a reduction of $\Bairespace$ to a closed set, $\restriction{(\pi
  \composition \psi)}{\Bairespace}$ is continuous, and $\sequence
  {\bfd_n}[n \in \N] \in \Heven[\N]$, then $(\pi \composition \psi)(\lim_{n
  \goesto \infty} \bfd_n) \notin \closure{\set{(\pi \composition \psi)
  (\bfd_n)}[n \in \N]}$, so the continuity of $\psi$ ensures that
  $\sequence{\pair{\psi(\bfd_n)}{(\pi \composition \psi)(\bfd_n)}}[n \in
  \N] \in H_\pi$, thus $(\restriction{\psi}{\Bairespace}) \times (\restriction
  {(\pi \composition \psi)}{\Bairespace})$ is a homomorphism from
  $\Heven[\N]$ to $H_\pi$.
  
  Conversely, suppose that $\phi \from \Bairespace \to \graph{\pi}$ is a
  continuous homomorphism from $\Heven[\N]$ to $H_\pi$, and set
  $\phi_X = \projection[X] \composition \phi$ and $\phi_Y = \projection[Y]
  \composition \phi$. As the definition of $H_\pi$ ensures that $\phi_X$
  is a homomorphism from $\Rcon[\N]{\compactification{\N}}$ to $\Rcon
  {X}$, Proposition \ref{compactification:extension} yields a continuous
  extension $\extension{\phi}_X \from \modifiedpower{\N}{\N}
  {\compactification{\N}} \to X$.

  If there exists $t \in \Bairetree$ for which $\restriction{\phi_Y}
  {\extensions{t}}$ is constant, then the function $\psi \from
  \modifiedpower{\N}{\N}{\compactification{\N}} \to X$ given by $\psi
  (\bfd) = \extension{\phi}_X(t \concatenation \bfd)$ is continuous and
  $\pi \composition \psi$ is a reduction of $\Bairespace$ to a singleton,
  so we can assume that
  \begin{equation} \tag{$\dagger$}
    \forall t \in \Bairetree \ \image{\phi_Y}{\extensions{t}} \text{ is infinite.}
  \end{equation}
  
  If there exists $t \in \Bairetree$ for which $\image{(\pi \composition
  \extension{\phi}_X)}{\extensions{t}[\compactification{\N}] \setminus
  \Bairespace}$ is finite, then by appealing to $(\dagger)$ and extending
  $t$, we can ensure that $\closure{\image{\phi_Y}{\extensions{t}}}$ and
  $\image{(\pi \composition \extension{\phi}_X)}{\extensions{t}
  [\compactification{\N}] \setminus \Bairespace}$ are disjoint, so the
  function $\psi \from \modifiedpower{\N}{\N}{\compactification{\N}} \to
  X$ given by $\psi(\bfd) = \extension{\phi}_X(t \concatenation \bfd)$ is
  continuous, $\pi \composition \psi$ is a reduction of $\Bairespace$ to
  $\closure{\image{\phi_Y}{\extensions{t}}}$, and $(\pi \composition \psi)
  (\bfd) = \phi_Y(t \concatenation \bfd)$ for all $\bfd \in \Bairespace$.
  We can therefore additionally assume that 
  \begin{equation} \tag{$\ddagger$}
    \forall t \in \Bairetree \ \image{(\pi \composition \extension{\phi}_X)}
      {\extensions{t}[\compactification{\N}] \setminus \Bairespace} \text
        { is infinite.}
  \end{equation}

  Fix positive real numbers $\epsilon_n \goesto 0$, as well as an
  enumeration $\sequence{t_k}[k \in \N]$ of $\Bairetree$ such that $t_j
  \extendedby t_k \implies j \le k$ for all $j, k \in \N$. We say that a
  sequence $\sequence{X_n}[n \in \N]$ of subsets of $X$ is \definedterm
  {closed and discrete} if there are no convergent sequences in
  $\product[n \in \N][X_n]$. We will recursively construct functions $f
  \from \Bairetree \to \Bairetree$ and $g \from \Bairetree \setminus \set
  {\emptysequence} \to \Bairetree$ such that, for all $k \in \N$, the
  following conditions hold:

  \begin{enumerate}
    \renewcommand{\theenumi}{\alph{enumi}}
    \item $\forall n \in \N \ f(t_k) \strictlyextendedby g(t_k \concatenation
      \sequence{n})$.
    \item $\forall m \neq n \ g(t_k \concatenation \sequence{m})(\length
      {f(t_k)}) \neq g(t_k \concatenation \sequence{n})(\length{f(t_k)})$.
    \item $\sequence{\image{\phi_Y}{\extensions{g(t_k \concatenation
      \sequence{n})}}}[n \in \N]$ is either closed-and-discrete or
        convergent.
    \item $g(t_{k+1}) \extendedby f(t_{k+1})$.
    \item $\forall j \le k \ (\pi \composition \extension{\phi}_X)(f(t_j)
      \concatenation \sequence{\infty}) \notin \closure{\image{\phi_Y}
        {\extensions{f(t_{k+1})}}}$.
    \item $\forall j \le k \ (\pi \composition \extension{\phi}_X)(f(t_{k+1})
      \concatenation \sequence{\infty}) \neq \lim_{n \goesto \infty} \image
        {\phi_Y}{\extensions{g(t_j \concatenation \sequence{n})}}$.
  \end{enumerate}
  We begin by setting $f(\emptysequence) = \emptysequence$.
  
  Suppose now that $k \in \N$ and we have already found $f(t_k)$, as
  well as $g(t_j \concatenation \sequence{n})$ for all $j < k$ and $n \in
  \N$, and fix extensions $\bfd_{k,n} \in \Bairespace$ of $f(t_k)
  \concatenation \sequence{n}$ for all $n \in \N$. Then there is an
  injective sequence $\sequence{i_{k,n}}[n \in \N]$ of natural numbers
  such that $\sequence{\set{\phi_Y(\bfd_{k,i_{k,n}})}}[n \in \N]$ is either
  closed-and-discrete or convergent, and continuity yields extensions
  $g(t_k \concatenation \sequence{n}) \strictlyextendedby \bfd_{k,
  i_{k,n}}$ of $f(t_k) \concatenation \sequence{i_{k,n}}$ with the property
  that $\image{\phi_Y}{\extensions{g(t_k \concatenation \sequence{n})}}
  \subseteq \ball[Y]{\phi_Y(\bfd_{k, i_{k,n}})}{\epsilon_n}$ for all $n \in \N$.
  By $(\dagger)$, there is an extension $\bfd_k \in \Bairespace$ of $g
  (t_{k+1})$ such that $\phi_Y(\bfd_k) \neq (\pi \composition \extension
  {\phi}_X)(f(t_j) \concatenation \sequence{\infty})$ for all $j \le k$, and
  therefore an extension $u_k \strictlyextendedby \bfd_k$ of $g(t_{k+1})$
  such that $(\pi \composition \extension{\phi}_X)(f(t_j) \concatenation
  \sequence{\infty}) \notin \closure{\image{\phi_Y}{\extensions{u_k}}}$
  for all $j \le k$. By $(\ddagger)$, there is an extension $f(t_{k+1})$ of
  $u_k$ satisfying condition (f), thereby completing the recursive
  construction.
  
  Conditions (a), (b), and (d) ensure that the function $\phi_X
  \composition \phi'$ is a homomorphism from $\Rcon[\N]
  {\compactification{\N}}$ to $\Rcon{X}$, where $\phi' \from \Bairespace
  \to \Bairespace$ is given by $\phi'(\bfd) = \union[n \in \N][f(\restriction
  {\bfd}{n})]$. Proposition \ref{compactification:extension} therefore
  yields a continuous extension $\psi \from \modifiedpower{\N}{\N}
  {\compactification{\N}} \to X$. As $\restriction{(\pi \composition \psi)}
  {\Bairespace} = \phi_Y \composition \phi'$, it only remains to show that
  $\pi \composition \psi$ is a reduction of $\Bairespace$ to a closed set.

  Suppose, towards a contradiction, that there exists $k \in \N$ with the
  property that $(\pi \composition \psi)(t_k \concatenation \sequence
  {\infty}) \in \closure{\image{(\pi \composition \psi)}{\Bairespace}}$.
  As conditions (a), (b), and (d) ensure that
  \begin{align*}
    (\pi \composition \psi)(t_k \concatenation \sequence{\infty})
      & \textstyle = \pi(\lim_{n \goesto \infty} \image{(\phi_X \composition
        \phi')} {\extensions{t_k \concatenation \sequence{n}}}) \\
      & \textstyle = \pi(\lim_{n \goesto \infty} \image{\phi_X}{\extensions
        {f(t_k \concatenation \sequence{n})}}) \\
      & = (\pi \composition \extension{\phi}_X)(f(t_k) \concatenation
        \sequence{\infty}),
  \end{align*}
  it follows that $\image{(\pi \composition \extension{\phi}_X)}{f(t_k)
  \concatenation \sequence{\infty}} \in \closure{\image{\phi_Y}
  {\Bairespace}}$. As condition (e) ensures that $\image{(\pi
  \composition \extension{\phi}_X)}{f(t_k) \concatenation \sequence
  {\infty}} \notin \closure{\image{\phi_Y}{\extensions{f(t_j)}}}$ for all $j >
  k$, and condition (f) implies that if $i < k$ and $\image{(\pi
  \composition \extension{\phi}_X)}{f(t_k) \concatenation \sequence
  {\infty}} \in \closure{\image{\phi_Y}{\extensions{f(t_i)}}}$ then there
  exists $j > i$ for which $\image{(\pi \composition \extension{\phi}_X)}
  {f(t_k) \concatenation \sequence{\infty}} \in \closure{\image{\phi_Y}
  {\extensions{f(t_j)}}}$, it follows that $\image{(\pi \composition
  \extension{\phi}_X)}{f(t_k) \concatenation \sequence{\infty}} \in
  \closure{\image{\phi_Y}{\extensions{f(t_k)}}} \setminus \union[n \in \N]
  [\closure{\image{\phi_Y}{\extensions{f(t_k \concatenation \sequence
  {n})}}}]$. Conditions (a), (b), and (d) therefore yield extensions $\bfc_n
  \in \Bairespace$ of $f(t_k) \concatenation \sequence{n}$ such that
  $\pi(\lim_{n \goesto \infty} \phi_X(\bfc_n)) \in \closure{\set{(\pi
  \composition \phi_X)(\bfc_n)}}$, contradicting the fact that $\phi$ is a
  homomorphism from $\Heven[\N]$ to $H_\pi$.
\end{propositionproof}

In particular, we obtain the following.

\begin{theorem}[{$\ogd[\N]{\Gamma}$}]
  \label{sigmacontinuous:dichotomy}
  Suppose that $X$ and $Y$ are separable metric spaces, $X \in
  \Gamma$, and $\pi \from X \to Y$ is \Borel. Then exactly one of the
  following holds:
  \begin{enumerate}
    \item The function $\pi$ is $\sigma$-continuous with closed
      witnesses.
    \item There is a continuous function $\psi \from \modifiedpower
      {\N}{\N}{\compactification{\N}} \to X$ with the property that $\pi
      \composition \psi$ is a reduction of $\Bairespace$ to a closed set
      and $\restriction{(\pi \composition \psi)}{\Bairespace}$ is continuous.
  \end{enumerate}
\end{theorem}

\begin{theoremproof}
  This follows from Proposition \ref{sigmacontinuous:hypergraph}.
\end{theoremproof}

The special case of Theorem \ref{sigmacontinuous:dichotomy}
where $\Gamma$ is the family of analytic spaces easily yields the
\Jayne-\Rogers characterization of the circumstances under which 
a function from an analytic metric space to a separable metric space
is $\sigma$-continuous with closed witnesses. The special case where
$\Gamma$ is the family of separable spaces easily yields the
generalization in which analyticity is weakened to separability under
the axiom of determinacy.

\section{Cardinal invariants and chromatic numbers} \label{invariants}

We begin this section by noting a straightforward restriction on chromatic
numbers of box-open hereditary $\N$-dimensional dihypergraphs.

\begin{proposition} \label{invariants:dominatingnumber}
  The chromatic number of $\Hevenprime[\N]$ is $\dominatingnumber$.
\end{proposition}

\begin{propositionproof}
  This follows from the fact that a subset of $\Bairespace$ is
  $\Hevenprime[\N]$-independent if and only if its closure is compact.
\end{propositionproof}

\begin{proposition}[{$\ogd[\N]{\Gamma}$}] \label{invariants:dichotomy}
  Suppose that $X$ is a \Hausdorff space in $\Gamma$ and $H$ is a
  box-open hereditary $\N$-dimensional dihypergraph on $X$. Then
  either $\chromaticnumber{H} \le \aleph_0$ or $\chromaticnumber{H}
  \ge \dominatingnumber$.
\end{proposition}

\begin{propositionproof}
  This follows from Proposition \ref{invariants:dominatingnumber}
  and the fact that every homomorphism from $\Heven[\N]$ to $H$
  is a homomorphism from $\Hevenprime[\N]$ to $H$.
\end{propositionproof}

As every analytic subset of a topological space is the union of
$\dominatingnumber$-many compact sets, the special case of
Proposition \ref{invariants:dichotomy} where $\Gamma$ is the family of
analytic sets and $H$ is of the form $H_X$ ensures that if an analytic
subset of a metric space is not contained in a \Ksigma set, then
$\dominatingnumber$ is the least cardinal $\kappa$ for which it is
contained in a union of $\kappa$-many compact sets.

As every analytic subset of a topological space is the union of
$\dominatingnumber$-many closed sets, the special case of Proposition
\ref{invariants:dichotomy} where $\Gamma$ is the family of analytic sets,
$H$ is of the form $H_{R,S}$, and $D = 1$ ensures that if an analytic
subset of a metric space cannot be separated from another set by an
\Fsigma set, then $\dominatingnumber$ is the least cardinal $\kappa$
for which it can be separated from the other set by a union of
$\kappa$-many closed sets.

As every analytic set that is the graph of a function is the union of
$\dominatingnumber$-many compact sets that are graphs of functions,
the special case of Proposition \ref{invariants:dichotomy} where
$\Gamma$ is the family of analytic sets and $H$ is of the form $H_\pi$
ensures that if a \Borel function $\pi$ from an analytic metric space to a
separable metric space is not $\sigma$-continuous with closed
witnesses, then $\dominatingnumber$ is the least cardinal $\kappa$ for
which $X$ is the union of $\kappa$-many closed sets on which $\pi$ is
continuous.

The special case of the last fact for \Baire-class-one functions is due to
\Solecki. He established this by noting that the two elements of his basis
for \Baire-class-one functions that are not $\sigma$-continuous have the
desired property. The above argument shows that this more
sophisticated basis theorem is unnecessary to obtain the desired result;
one need only observe that $\sigma$-continuity with closed witnesses
can be characterized using a box-open hereditary $\N$-dimensional
dihypergraph.

We next turn our attention to the computation of the chromatic number
of $\Heven[\N]$ itself. For each partial function $f \from \functions{<\N}
{D} \partialto D$, define $\bfD_f = \set{\bfd \in \functions{\N}{D}}[\forall
n \in \N \ (\restriction{\bfd}{n} \in \domain{f} \implies \bfd(n) \neq
f(\restriction{\bfd}{n}))]$.

\begin{proposition} \label{invariants:family}
  Suppose that $D$ is a set of cardinality at least two. Then
  $\chromaticnumber{\Heven[D]} = \min \set{\cardinality{\calF}}[{\calF
  \subseteq \functions{\functions{<\N}{D}}{D} \mathand \functions{\N}{D}
  = \union[f \in \calF][\bfD_f]}]$.
\end{proposition}

\begin{propositionproof}
  This follows from the fact that a set $X \subseteq \functions{\N}{D}$
is $\Heven[D]$-independent if and only if there is a function $f
\from \functions{<\N}{D} \to D$ for which $X \subseteq \bfD_f$.
\end{propositionproof}

We now establish an analog of Proposition \ref{invariants:dichotomy}
without the assumption that $H$ is hereditary.

\begin{proposition}[{$\ogd[D]{\Gamma}$}] \label{invariants:bounds}
  Suppose that $D$ is a countable set of cardinality at least two, $X$
  is a \Hausdorff space in $\Gamma$, and $H$ is a box-open
  $D$-dimensional dihypergraph on $X$. Then either $\chromaticnumber
  {H}\le \aleph_0$ or $\chromaticnumber{H} \ge \coveringnumber
  {\meagerideal}$.
\end{proposition}

\begin{propositionproof}
  It is sufficient to show that $\chromaticnumber{\Heven[D]} \ge
  \coveringnumber{\meagerideal}$, which follows from Proposition \ref
  {invariants:family} and the observation that if $f \from \functions{<\N}
  {D} \to D$, then $\bfD_f$ is meager with respect to the usual topology
  on $\functions{\N}{D}$.
\end{propositionproof}

Alternatively, one can obtain the above result by noting that if $H$ is a
box-open $D$-dimensional dihypergraph on a \Hausdorff space,
$\calI_H$ is the $\sigma$-ideal generated by the family of closed
$H$-independent sets, and $\chromaticnumber{H} > \aleph_0$, then
$\chromaticnumber{H} = \coveringnumber{\calI_H}$.

Conversely, given a $\sigma$-ideal $\calI$ on a topological space $X$,
let $H_\calI$ be the $\N$-dimensional dihypergraph on $X$ consisting
of all sequences $\sequence{x_n}[n \in \N]$ of elements of $X$ for
which $\closure{\set{x_n}[n \in \N]} \notin \calI$.

\begin{proposition} \label{invariants:hypergraph}
  Suppose that $X$ is a separable metric space and $\calI$ is a
  $\sigma$-ideal on $X$ generated by a family of closed subsets of
  $X$.
  \begin{enumerate}
    \item If $\calI$ covers $X$, then $H_\calI$ is box open.
    \item If $\coveringnumber{\calI} > \aleph_0$, then
      $\chromaticnumber{H_\calI} = \coveringnumber{\calI}$.
  \end{enumerate}
\end{proposition}

\begin{propositionproof}
  To see (1), suppose that $\sequence{x_n}[n \in \N] \in H_\calI$, fix
  positive real numbers $\epsilon_n \goesto 0$, and observe that if
  $\sequence{y_n}[n \in \N] \in \product[n \in \N][{\ball[X]{x_n}{\epsilon_n}}]$,
  then $\closure{\set{x_n}[n \in \N]} \subseteq \set{x_n}[n \in \N] \union
  \closure{\set{y_n}[n \in \N]}$, thus $\sequence{y_n}[n \in \N] \in H_\calI$.
  
  To see (2), note first that every $H_\calI$-independent set is in $\calI$,
  so $\coveringnumber{\calI} \le \chromaticnumber{H_\calI}$.
  Conversely, as every set in $\calI$ is contained in the union of
  countably-many closed sets in $\calI$, the fact that every closed set in
  $\calI$ is $H_\calI$-independent ensures that $\chromaticnumber
  {H_{\calI}} \le \coveringnumber{\calI} \cdot \aleph_0 = \coveringnumber
  {\calI}$. 
\end{propositionproof}

As a corollary, we also obtain an analog of Proposition \ref
{invariants:dominatingnumber} without the assumption that $H$ is
hereditary.

\begin{proposition}
  The chromatic number of $\Heven[\N]$ is $\coveringnumber
  {\meagerideal}$.
\end{proposition}

\begin{propositionproof}
  Proposition \ref{invariants:bounds} ensures that $\chromaticnumber
  {\Heven[\N]} \ge \coveringnumber{\meagerideal}$. As the open
  dihypergraph dichotomy yields a homomorphism from $\Heven[\N]$ to
  $H_\calM$, Proposition \ref{invariants:hypergraph} implies that
  $\chromaticnumber{\Heven[\N]} \le \coveringnumber{\meagerideal}$.
\end{propositionproof}

Stronger bounds can be obtained when $D$ is finite. In the special
case that $D = 2$, this is trivial.

\begin{proposition}[{$\ogd{\Gamma}$}]
  \label{invariants:continuum:ogd}
  Suppose that $X$ is a \Hausdorff space in $\Gamma$ and $G$ is an
  open graph on $X$. Then either $\chromaticnumber{G} \le \aleph_0$
  or $\chromaticnumber{G} = \continuum$.
\end{proposition}

\begin{propositionproof}
  It is sufficient to show that $\chromaticnumber{\Heven[2]} =
  \continuum$, which follows from the fact that $\Heven[2] =
  \completegraph{\Cantorspace}$.
\end{propositionproof}

The case that $2 < \cardinality{D} < \aleph_0$ is substantially subtler,
and related to the adaptive global prediction numbers considered in
\cite[\S10]{Blass}. Our first observation is a straightforward analog of
Proposition \ref{invariants:bounds}.

\begin{proposition}[{$\ogd[D]{\Gamma}$}]
  \label{invariants:covN:ogd}
  Suppose that $D$ is a finite set of cardinality at least two, $X$ is a
  \Hausdorff space in $\Gamma$, and $H$ is an open $D$-dimensional
  dihypergraph on $X$. Then either $\chromaticnumber{H} \le \aleph_0$
  or $\chromaticnumber{H} \ge \coveringnumber{\nullideal}$.
\end{proposition}

\begin{propositionproof}
  It is sufficient to show that $\chromaticnumber{\Heven[D]} \ge
  \coveringnumber{\nullideal}$, which follows from Proposition \ref
  {invariants:family} and the observation that if $f \from \functions{<\N}
  {D} \to D$ and $\mu$ is a strictly positive probability measure on
  $D$, then $\bfD_f$ is $\functions{\N}{\mu}$-null.
\end{propositionproof}

We next establish an analog of the main result of \cite{Brendle}. We
say that a partial function $u \from \N \partialto D$ is a \definedterm
{$D$-\Silver condition} if its domain is co-infinite, and we say that a set
$X \subseteq \functions{\N}{D}$ is \definedterm{$D$-\Silver null} if
every $D$-\Silver condition $u \from \N \partialto D$ extends to a
$D$-\Silver condition $v \from \N \partialto D$ for which $\extensions
{v}$ is disjoint from $X$.

\begin{theorem} \label{invariants:boundingnumber:Heven}
  Suppose that $D$ is a finite set of cardinality at least two, $X
  \subseteq \functions{\N}{D}$, and $\chromaticnumber{\restriction
  {\Heven[D]}{X}} < \boundingnumber$. Then $X$ is $D$-\Silver null.
\end{theorem}

\begin{theoremproof}
  By Proposition \ref{invariants:family}, it is sufficient to show that if
  $\calF \subseteq \functions{\functions{<\N}{D}}{D}$ has cardinality
  strictly less than $\boundingnumber$, then $\union[f \in \calF][\bfD_f]$
  is $D$-\Silver null. By a straightforward recursive construction of
  length $\cardinality{D}$, we need only show that for all $d \in D$,
  every $D$-\Silver condition $u \from \N \partialto D$ extends to a
  $D$-\Silver condition $v \from \N \partialto D$ with the property that
  \begin{equation*}
    \forall f \in \calF \forall \bfd \in \bfD_f \intersection \extensions{v}
      \forcofinitelymany n \in \setcomplement{\domain{v}} \ d \neq
        f(\restriction{\bfd}{n}).
  \end{equation*}
  It is clearly sufficient to handle the special case that $u =
  \emptysequence$.
  
  For all $f \in \calF$, fix a function $g_f \from \N \to \N$ such that for all
  $n \in \N$ and $t \in \functions{n}{D}$, if there exists $i \in \N$ for which
  $d = f(t \concatenation \constantsequence{d}{i})$, then there is such
  an $i < g_f(n)$. Fix a function $g \from \N \to \N$ eventually dominating
  $g_f$ for all $f \in \calF$, set $h'(0) = 0$, and recursively define $h(n) =
  g(h'(n) + 1)$ and $h'(n + 1) = h'(n) + 1 + h(n)$ for all $n \in \N$. Let $v$
  be the function on $\setcomplement{\set{h'(n+1)}[n \in \N]}$ with
  constant value $d$.
  
  To see that $v$ is as desired, note first that if $\bfd \in \extensions{v}$,
  then a straightforward inductive argument reveals that $\restriction
  {\bfd}{h'(n+1)}$ is the concatenation of $\restriction{\bfd}{(h'(n) + 1)}$
  and $\constantsequence{d}{h(n)}$, for all $n \in \N$. It follows that if
  $f \in \calF$, $\bfd \in \bfD_f$, and $n \in \N$ is sufficiently large that
  $g_f(h'(n) + 1) \le g(h'(n) + 1)$, then $d \neq f(\restriction{\bfd}
  {h'(n+1)})$.
\end{theoremproof}

\begin{theorem}[{$\ogd[D]{\Gamma}$}]
  \label{invariants:boundingnumber:ogd}
  Suppose that $D$ is a finite set of cardinality at least two, $X$ is a
  \Hausdorff space in $\Gamma$, and $H$ is an open $D$-dimensional
  dihypergraph on $X$. Then either $\chromaticnumber{H} \le \aleph_0$
  or $\chromaticnumber{H} \ge \boundingnumber$.
\end{theorem}

\begin{theoremproof}
  It is sufficient to show that $\chromaticnumber{\Heven[D]} \ge
  \boundingnumber$, which follows from Theorem \ref
  {invariants:boundingnumber:Heven}.
\end{theoremproof}

We now establish an analog of the main result of \cite{Kamo}.

\begin{theorem} \label{invariants:cofinality:Heven}
  It is consistent that $\chromaticnumber{\Heven[D]} > \cofinality
  {\nullideal}$ for every finite set $D$ of cardinality at least two.
\end{theorem}

\begin{theoremproof}
  Given a natural number $D \ge 2$, we say that a subtree $T$ of
  $\functions{<\N}{D}$ is \definedterm{$D$-perfect} if every element of
  $T$ has an extension $t \in T$ with the property that $\forall d \in D \ t
  \concatenation \sequence{d} \in T$. Let $S_D$ denote the set of
  such trees. Let $\le$ denote the partial order on $S_D$ with respect to
  which $S \le T$ if and only if $S \subseteq T$, and define $\bbS_D =
  \pair{S_D}{\le}$. When $D = 2$, this is just \Sacks forcing, and the usual
  proof that the latter is proper and has the \Sacks property works just
  as well for every $\bbS_D$.
  
  \begin{lemma} \label{invariant:cofinality:Heven:Sacks}
    Suppose that $D \ge 2$ is a natural number, $G$ is
    $\bbS_D$-generic over $V$, $\bfd$ is the unique branch through
    $\intersection[G]$, and $f \from \functions{<\N}{D} \to D$ is in $V$.
    Then there exists $n \in \N$ for which $\bfd(n) = f(\restriction{\bfd}
    {n})$.
  \end{lemma}
  
  \begin{lemmaproof}
    This follows from the fact that there are $\le$-densely many trees
    $T \in S_D$ for which there exists $t \in \functions{<\N}{D}$ such
    that $t \concatenation \sequence{f(t)}$ is $\extendedby$-comparable
    with every element of $T$.
  \end{lemmaproof}
 
  Suppose now that $\CH$ holds in $V$, and fix a sequence $\sequence
  {D_\beta}[\beta < \omega_2]$ of natural numbers that are at least two
  for which every such natural number appears cofinally often, as well
  as a countable support iteration $\sequence{\bbP_\alpha, \dot
  {\bbQ}_\beta}[\alpha \le \omega_2, \beta < \omega_2]$ such that
  $\forces[\bbP_\beta] \dot{\bbQ}_\beta = (\bbS_{D_\beta})^
  {V^{\bbP_\beta}}$ for all $\beta < \omega_2$.
 
  As $\bbP_{\omega_2}$ has the Sacks property, it follows that
  $\cofinality{\nullideal} = \aleph_1$ in $V^{\bbP_{\omega_2}}$ (see, for
  example, \cite[\S11.5]{Blass}). To see that $\chromaticnumber{\Heven
  [D]} > \aleph_1$ in $V^{\bbP_{\omega_2}}$ for all natural numbers $D
  \ge 2$, note that if $\calF \subseteq \functions{\functions{<\N}{D}}{D}$
  is in $V^{\bbP_{\omega_2}}$ and has cardinality at most $\aleph_1$ in
  $V^{\bbP_{\omega_2}}$, then there exists $\alpha < \omega_2$ for
  which $\calF \in V^{\bbP_\alpha}$. Fix $\beta \ge \alpha$ for which
  $D_\beta = D$, note that $\functions{\N}{D} \nsubseteq \union[f \in
  \calF][\bfD_f]$ in $V^{\bbP_{\beta+1}}$ by Lemma \ref
  {invariant:cofinality:Heven:Sacks}, and appeal to Proposition \ref
  {invariants:family}.
\end{theoremproof}

\begin{theorem}[{$\ogd[D]{\Gamma}$}]
  \label{invariants:cofinality:ogd}
  It is consistent that whenever $D$ is a finite set of cardinality at least
  two, $X$ is a \Hausdorff space in $\Gamma$, and $H$ is an open
  $D$-dimensional dihypergraph on $X$, either $\chromaticnumber{H}
  \le \aleph_0$ or $\chromaticnumber{H} > \cofinality{\nullideal}$.
\end{theorem}

\begin{theoremproof}
  It is sufficient to establish the consistency of $\chromaticnumber
  {\Heven[D]} > \cofinality{\nullideal}$ for all finite sets $D$ of cardinality
  at least two, which follows from Theorem \ref
  {invariants:cofinality:Heven}.
\end{theoremproof}

We finally establish a consistent upper bound as well.

\begin{theorem} \label{invariants:dominatingnumber:strictlybelow}
  It is consistent that $\chromaticnumber{\Heven[D]} <
  \dominatingnumber$ for every finite set $D$ of cardinality at least
  three.
\end{theorem}

\begin{theoremproof}
  Given a natural number $D \ge 3$, let $P_D$ denote the set of pairs
  $\pair{\bfD}{f}$ with the property that $f \from \functions{<n}{D} \to D$
  for some $n \in \N$, $\bfD \subseteq \bfD_f$, and $\restriction{\bfc}{n}
  \neq \restriction{\bfd}{n}$ for all distinct $\bfc, \bfd \in \bfD$. Let $\le$
  denote the partial order on $P_D$ with respect to which $\pair{\bfC}{f}
  \le \pair{\bfD}{g}$ if and only if $\bfC \subseteq \bfD$ and $f
  \extendedby g$, and define $\bbP_D = \pair{P_D}{\le}$.
  
  \begin{lemma} \label{invariants:dominatingnumber:strictlybelow:linked}
    Suppose that $n \in \N$ and $f \from \functions{<n}{D} \to D$. Then
    any two elements of $P_D$ of the form $\pair{\bfC}{f}$ and $\pair{\bfD}
    {f}$ are $\le$-compatible.
  \end{lemma}
  
  \begin{lemmaproof}
    Fix $m \ge n$ sufficiently large that $\restriction{c}{m} \neq \restriction
    {d}{m}$ for all distinct $c \in \bfC$ and $d \in \bfD$. As $D \ge 3$,
    there exists $g \from \functions{<m}{D} \to D$ such that $\bfC \union
    \bfD \subseteq \bfD_g$ and $f \extendedby g$, in which case $\pair
    {\bfC \union \bfD}{g}$ is a common $\le$-extension of $\pair{\bfC}{f}$
    and $\pair{\bfD}{f}$.
  \end{lemmaproof}
  
  Fix a sequence $\sequence{D_n}[n \in \N]$ of natural numbers that
  are at least three for which every such natural number appears
  cofinally often, and let $\bbP = \pair{P}{\le}$ denote the finite
  support product of $\sequence{\bbP_{D_n}}[n \in \N]$. Lemma \ref
  {invariants:dominatingnumber:strictlybelow:linked} ensures that
  each of the partial orders $\bbP_D$ is $\sigma$-linked, thus so too
  is $\bbP$.
  
  \begin{lemma}
    \label{invariants:domintatingnumber:strictlybelow:coloring}
    Suppose that $D \ge 3$ is a natural number, $G$ is $\bbP$-generic
    over $V$, and $\bfd \in (\functions{\N}{D})^V$. Then there exists $k \in
    \N$ such that $D_k = D$ and $\bfd(n) \neq f(\restriction{\bfd}{n})$ for
    all $n \in \N$, where $f = \union[{\set{f_k}[{\sequence{\bfD_n, f_n}[n
    \in \N] \in G}]}]$.
  \end{lemma}
  
  \begin{lemmaproof}
    This follows from the fact that there are $\le$-densely many
    sequences $\sequence{\bfD_n, f_n}[n \in \N] \in P$ for which
    there exists $k \in \N$ such that $D = D_k$ and $\bfd \in \bfD_k$.
  \end{lemmaproof}
  
  \begin{lemma} \label{invariants:dominatingnumber:strictlybelow:bound}
    Suppose that $\dot{\bfd}$ is a $\bbP$-name for an element of
    $\Bairespace$, $\sequence{k_n}[n \in \N] \in \Bairespace$, and
    $\sequence{f_n}[n \in \N]$ is in the finite support product of
    $\sequence{\functions{\functions{<k_n}{D_n}}{D_n}}[n \in \N]$. Then
    $\forall i \in \N \exists j \in \N \forall \sequence{\bfD_n, f_n}[n \in \N] \in
    P \ \sequence{\bfD_n, f_n}[n \in \N] \nforces[\bbP] \dot{\bfd}(i) \ge j$.
  \end{lemma}
  
  \begin{lemmaproof}
    Suppose, towards a contradiction, that there exist $i \in \N$ and
    $\sequence{\bfD_{j,n}, f_n}[n \in \N] \in P$ such that $\sequence
    {\bfD_{j,n}, f_n}[n \in \N] \forces[\bbP] \dot{\bfd}(i) \ge j$ for all $j \in
    \N$. By passing to a subsequence, we can assume that there
    are sequences $\sequence{I_n}[n \in \N]$ of finite sets and
    $\sequence{\bfd_{i,j,n}}[\triple{i}{j}{n} \in I_n \times \N \times \N]$ of
    elements of $\Bairespace$ such that $\bfD_{j, n} = \set{\bfd_{i,j,n}}[i
    \in I_n]$ for all $j, n \in \N$ and $\bfd_{i,n} = \lim_{j \goesto \infty}
    \bfd_{i,j,n}$ exists for all $i, n \in \N$. Set $\bfD_n = \set{\bfd_{i,n}}
    [i \in I_n]$. Then $\sequence{\bfD_n, f_n}[n \in \N] \in P$, so there
    exist $k \in \N$ and an extension of $\sequence{\bfD_n, f_n}[n \in
    \N] \in P$ for which $\sequence{\bfD_n', f_n'}[n \in \N]
    \forces[\bbP] \dot{\bfd}(i) = k$. As $\sequence{\bfD_{j,n}, f_n'}[n \in
    \N] \in P$ for all sufficiently large $j \in \N$, Lemma \ref
    {invariants:dominatingnumber:strictlybelow:linked} ensures that
    $\sequence{\bfD_{j,n}, f_n}[n \in \N]$ and $\sequence{\bfD_n', f_n'}
    [n \in \N]$ are compatible for all sufficiently large $j > k$, the
    desired contradiction.
  \end{lemmaproof}
  
  Suppose now that $\dominatingnumber > \aleph_1$ in $V$, and fix
  a finite support iteration $\sequence{\bbP_\alpha, \dot{\bbQ}_\beta}
  [\alpha \le \omega_1, \beta < \omega_1]$ such that $\forces
  [\bbP_\beta] \dot{\bbQ}_\beta = \bbP^{V^{\bbP_\beta}}$ for all
  $\beta < \omega_1$.
  
  Proposition \ref{invariants:family} and Lemma \ref
  {invariants:domintatingnumber:strictlybelow:coloring} easily imply
  that $\chromaticnumber{\Heven[D]} = \aleph_1$ in
  $V^{\bbP_{\omega_1}}$. To see that $\dominatingnumber >
  \aleph_1$ in $V^{\bbP_{\omega_1}}$, note that if $\sequence{\dot
  {\bfd}_\alpha}[\alpha < \omega_1] \in V$ is a sequence of
  $\bbP_{\omega_1}$-names for elements of $\Bairespace$, then
  Lemma \ref{invariants:dominatingnumber:strictlybelow:bound}
  yields a sequence $\sequence{\bfd_\alpha}[\alpha < \omega_1] \in
  V$ of elements of $\Bairespace$ such that $\forces
  [\bbP_{\omega_1}] \dot{\bfd}_\alpha \eventuallydominatedby
  \bfd_\alpha$ for all $\alpha < \omega_1$. The fact that
  $\dominatingnumber > \aleph_1$ in $V$ then yields $\bfd \in
  (\Bairespace)^V$ such that $\bfd \neventuallydominatedby
  \bfd_\alpha$ for all $\alpha < \omega_1$, so $\forces[\bbP_
  {\omega_1}] \bfd \neventuallydominatedby \dot{\bfd}_\alpha$ for
  all $\alpha < \omega_1$.
\end{theoremproof}

\providecommand{\bysame}{\leavevmode\hbox to3em{\hrulefill}\thinspace}
\providecommand{\MR}{\relax\ifhmode\unskip\space\fi MR }
\providecommand{\MRhref}[2]{%
  \href{http://www.ams.org/mathscinet-getitem?mr=#1}{#2}
}
\providecommand{\href}[2]{#2}

\end{document}